\documentclass[11pt]{amsart}

\usepackage{amsmath,latexsym}
\usepackage{amssymb}
\usepackage{color}

\definecolor{Myblue}{rgb}{0.0,0,0.9}
\definecolor{Mygreen}{rgb}{0.2,1,0}

\newcommand{\EF}[1]{{\color{cyan}#1}}
\newcommand{\MB}[1]{{\color{red}#1}}
\newcommand{\ST}[1]{{\color{green}#1}}
\newcommand{\IP}[1]{{\color{magenta}#1}}
\newcommand{\MIP}[1]{{\color{Myblue}#1}}
\newcommand{\MP}[1]{{\color{yellow}#1}}

\newtheorem{thm}{Theorem}[section]

\newtheorem{coro}[thm]{Corollary}

\theoremstyle{definition}
\newtheorem{defn}[thm]{Definition}
\theoremstyle{definition}
\newtheorem{exmp}[thm]{Example}

\numberwithin{equation}{section}

\newcommand{\La}{A}
\newcommand{\orho}{\overline{\rho}}
\newcommand{\Ga}{\Gamma}
\newcommand{\GaLa}{\Gamma_A}
\newcommand{\Gaa}{\Gamma$(mod$A)}
\newcommand{\bl}{\bf\huge}
\newcommand{\lla}{{\cal L}_{\La}}
\newcommand{\rla}{{\cal R}_{\La}}
\newcommand{\llb}{{\cal L}_{B}}
\newcommand{\rlb}{{\cal R}_{B}}
\newcommand{\TT}{{\cal T}(T)}
\newcommand{\FT}{{\cal F}(T)}

\newcommand{\ola}[1]{\overleftarrow{#1}}
\newcommand{\llbl}{{\cal L}_{B'}}
\newcommand{\rlbl}{{\cal R}_{B'}}
\newcommand{\vs}{\vspace{.3 cm}}
\newcommand{\vsd}{\vspace{.2 cm}}
\newcommand{\pn}{\par\noindent}
\newcommand{\lra}{\longrightarrow}
\newcommand{\st}{\stackrel}
\newcommand{\ipla}{${\mathcal P}^g_\La$}
\newcommand{\ipb}{${\mathcal P}^g_B$}
\newcommand{\iplan}{${\mathcal P}^g_\La \neq \emptyset$}
\newcommand{\ipbn}{${\mathcal P}^g_B \neq \emptyset$}
\newcommand{\benu}{\begin{enumerate}}
\newcommand{\enu}{\end{enumerate}}
\newcommand{\po}{P_{\omega}}
\newcommand{\sm}{{\cal S} (\rightarrow M)}
\newcommand{\bema}{\left[\begin{array}}
\newcommand{\ema}{\end{array}\right]}
\newcommand{\Al}{A_\lambda}
\newcommand{\Ar}{A_\rho}
\newcommand{\gen}[1]{\langle#1\rangle}
\newcommand{\field}{k}
\newcommand{\LL}{\begin{picture}(15,5)\put(2,3){\line(1,0){10}}\end{picture}}
\newcommand{\xLL}[1]{\begin{picture}(15,12)\put(2,3){\line(1,0){10}}\put(7.5,6){\HBCenter{\small
      $#1$}}\end{picture}}

\newcommand{\CluCat}{\mathcal{C}}

\renewcommand{\mod}{\operatorname{mod}}
\newcommand{\op}{\operatorname{op}}
\newcommand{\add}{\operatorname{add}}
\newcommand{\End}{\operatorname{End}}
\newcommand{\Ind}{\operatorname{Ind}}
\newcommand{\Hom}{\operatorname{Hom}}
\newcommand{\Ext}{\operatorname{Ext}}
\newcommand{\Ker}{\operatorname{Ker}}
\newcommand{\Der}{\operatorname{D^b}}
\newcommand{\HomD}{\operatorname{Hom}_{\Der(H)}}
\newcommand{\HomC}{\operatorname{Hom}_{\CluCat}}
\newcommand{\EndD}{\operatorname{End}_{\Der(H)}}
\newcommand{\EndC}{\operatorname{End}_{\CluCat}}
\newcommand{\dual}{\mathrm{D}}
\newcommand{\ealg}{\mathcal{R}}
\newcommand{\Hh}{\mathcal{H}}
\newcommand{\Pp}{\mathcal{P}}
\newcommand{\Rr}{\mathcal{R}}
\newcommand{\Cc}{\mathcal{C}}
\newcommand{\Tt}{\mathcal{T}}
\newcommand{\Ii}{\mathcal{I}}
\newcommand{\rad}{\operatorname{rad}}
\newcommand{\soc}{\operatorname{soc}}
\newcommand{\TOP}{\operatorname{top}}
\newcommand{\ind}{\operatorname{ind}}
\newcommand{\gldim}{\operatorname{gldim}}
\newcommand{\pdim}{\operatorname{pdim}}
\newcommand{\Groth}[1]{\operatorname{K_\circ}(#1)}
\newcommand{\proj}{\mathbf{P}}
\newcommand{\inj}{\mathbf{I}}
\newcommand{\An}{\mathbb{A}}
\newcommand{\ZZ}{\mathbb{Z}}
\newcommand{\ra}{\rightarrow}
\newcommand{\oT}{\widetilde{T}}
\newcommand{\oB}{\widetilde{B}}
\newcommand{\id}{\operatorname{id}}
\newcommand{\trivpath}[1]{\operatorname{e}_{#1}}

\def\semidirprod{\begin{picture}(8,8)\qbezier(2,0.5)(5,3.5)(8,6.5)\qbezier(2,6.5)(5,3.5)(8,0.5)\put(2,0.5){\line(0,1){6}}\end{picture}}
\def\ssemidirprod{{\unitlength=0.8pt\semidirprod}}
\def\sssemidirprod{{\unitlength=0.7pt\semidirprod}}

\newcommand{\proofend}{\hfill$\Box$\par}

\newcommand{\HVCenter}[1]{\setbox 0=\hbox{#1}%
        \dimen0=\wd0%
        \dimen1=\ht0%
        \divide\dimen0 by 2%
        \divide\dimen1 by 2%
        \hskip -\dimen0%
        \lower \dimen1%
        \box0%
        \hskip -\dimen0}
\newcommand{\HBCenter}[1]{\setbox 0=\hbox{#1}%
        \dimen0=\wd0%
        \dimen1=\ht0%
        \divide\dimen0 by 2%
        \hskip -\dimen0%
        \box0%
        \hskip -\dimen0}
\newcommand{\LTCenter}[1]{\setbox 0=\hbox{#1}%
        \dimen1=\ht0%
        \lower \dimen1%
        \box0%
        \hskip -\dimen0}
\newcommand{\HTCenter}[1]{\setbox 0=\hbox{#1}%
        \dimen0=\wd0%
        \dimen1=\ht0%
        \divide\dimen0 by 2%
        \hskip -\dimen0%
        \lower \dimen1%
        \box0%
        \hskip -\dimen0}
\newcommand{\RTCenter}[1]{\setbox 0=\hbox{#1}%
        \dimen0=\wd0%
        \dimen1=\ht0%
        \hskip -\dimen0%
        \lower \dimen1%
        \box0%
        \hskip -\dimen0}
\newcommand{\RBCenter}[1]{\setbox 0=\hbox{#1}%
        \dimen0=\wd0%
        \dimen1=\ht0%
        \hskip -\dimen0%
        \box0%
        \hskip -\dimen0}
\newcommand{\RVCenter}[1]{\setbox 0=\hbox{#1}%
        \dimen0=\wd0%
        \dimen1=\ht0%
        \divide\dimen1 by 2%
        \hskip -\dimen0%
        \lower \dimen1%
        \box0%
        \hskip -\dimen0}
\newcommand{\LVCenter}[1]{\setbox 0=\hbox{#1}%
        \dimen1=\ht0%
        \divide\dimen1 by 2%
        \lower \dimen1%
        \box0%
        \hskip -\dimen0}

\newcommand{\mylabel}[1]{\put(0,0){\color{white} \circle*{7}}\put(0,0){\circle{7}}\put(0,0){\HVCenter{\tiny\bf #1}}}
\newcommand{\mylabeltwo}[2]{\put(0,0){\color{white} \circle*{7}}\put(0,0){\circle{7}}\put(0,0){\HVCenter{\tiny\bf #1}}\put(0,4.5){\HVCenter{\tiny #2}}}

\newcommand{\mypdot}{\put(0,0){\color{white} \circle*{4}}\put(0,0){\circle*{2}}}
\newcommand{\mypdott}{\put(0,0){\circle*{2}}}

\begin{document}
\sloppy

\title[Quasitilted algebras  \, ]{
\small {\sc On the quiver with relations of a quasitilted algebra
and applications }}

\author[Bordino]{Natalia Bordino}
\address{Natalia Bordino,
Departamento de Matem\'atica, Facultad de Ciencias Exactas y
  Naturales, Funes 3350, Universidad Nacional de Mar del Plata, 7600
  Mar del Plata, Argentina.}
\email{bordino@mdp.edu.ar}

\author[Fern\'andez]{Elsa Fern\'andez}
\address{Elsa Fern\'andez,
Facultad de Ingenier\'{\i}a, Universidad Nacional de la
  Patagonia San Juan Bosco, 9120 Puerto Madryn, Argentina.}
\email{elsafer9@gmail.com}

\author[Trepode]{Sonia Trepode}
\address{Sonia Trepode,
Departamento de Matem\'atica, Facultad de Ciencias Exactas y
  Naturales, Funes 3350, Universidad Nacional de Mar del Plata, 7600
  Mar del Plata, Argentina.}
\email{strepode@mdp.edu.ar}

\begin{abstract}
In this paper we discuss, in terms of quivers with relations,
sufficient and necessary conditions for an algebra to be a
quasitilted algebra. We start with an algebra with global dimension
two and we give a sufficient condition for it to be a quasitilted
algebra. We show that this condition is not necessary. In the case
of a strongly simply connected schurian algebra, we discuss
necessary conditions, and, combining both type of conditions, we are
able to analyze if the given algebra is quasitilted. As an
application we obtain the quiver with relations of all the tilted
and cluster tilted algebras of Dynkin type $E_p$.
\end{abstract}

\thanks{We thank M. I. Platzeck for useful discussions about
Theorem \ref{thm:ct1}. This work is part of the Ph.D. Thesis of
Natalia Bordino, under the supervision of her advisers, Sonia
Trepode and Elsa Fern\'andez. The third author is a researcher of
CONICET, Argentina.}
\subjclass[2000]{Primary:
16G20, 
Secondary: 16G70, 
} \maketitle


\section{Introduction}

An interesting problem in the theory of representation of algebras
is to know if an algebra given by a quiver with relations is a
tilted algebra. This problem was studied by several authors and has
been solved in particular cases, see \cite{a}, \cite{cpt1},
\cite{cpt2}, \cite{cpt3}, \cite{h}, \cite{hl1}, \cite{hl2},
\cite{hr}. In this work, we consider a bigger class of algebras,
that is, the quasitilted algebras introduced by Happel, Reiten and
Smalo, see \cite{hrs}.

We introduced a notion of bounded consecutive relations inspired by the one
given by Assem and Redondo in \cite{ar}. Our first result states that, an algebra
with global dimension two having non consecutive bounded relations, is a
quasitilted algebra. Therefore, we show that the converse of this result does
not hold. On the other hand, if a quasitilted  algebra has bounded consecutive relations,
we give a necessary condition over this kind of relations.

Finally, we obtain a sufficient condition for a strongly simply
connected schurian algebra to be a quasitilted algebra. We arrive at
this result combining the criterion for global dimension two,
developed in \cite{bft} with the previous results.

As applications of our main result, and using the results given in \cite{f}, we
are able to give all the quivers with relations of the tilted and cluster tilted algebras of Dynkin type.
In this work we consider tilted an cluster tilted algebras of type $E_p$.

\section{Preliminaries}

In this paper, by an algebra, we always mean a basic and connected
finite dimensional algebra over an algebraically closed field $k$.
Given a quiver $Q$, we denote by $Q_0$ its set of vertices and by
$Q_1$ its set of arrows. A \textit{relation} in $Q$ from a vertex
$x$ to a vertex $y$ is a linear combination $\rho = \sum^m_{i=1}
\lambda_i w_i$ where, for each $i$, $\lambda_i \in k$ is non-zero
and $w_i$ is a path of length at least two from $x$ to $y$.

We denote by $kQ$ the path algebra of $Q$ and by $kQ(x, y)$ the
$k$-vector space generated by all paths in $Q$ from $x$ to $y$.
For an algebra $A$, we denote by $Q_A$ its ordinary quiver. For
every algebra $A$, there exists an ideal $I$ in $kQ_A$, generated
by a set of relations, such that $A \simeq kQ_A/I$. The pair
$(Q_A, I)$ is called a \textit{presentation} of $A$ and $A$ is
said to be given by the \textit{bound quiver} $(Q_A, I)$.

An algebra $A$ is called \textit{triangular} if $Q_A$ has no
oriented cycles, and it is called \textit{schurian} if, for all
$x, y \in A_0$, we have $\mbox{dim}_k A(x, y) \leq 1$. A
triangular algebra $A$ is called \textit{simply connected} if, for
any presentation $(Q_A, I)$ of $A$, the group $\pi_1(Q_A, I)$ is
trivial, see \cite{bg}. It is called \textit{strongly simply
connected} if every full convex subcategory of $A$ is simply
connected, \cite{Sko}.

In this work, we always deal with schurian triangular algebras. For
a vertex $x$ in the quiver $Q_A$, we denote by $e_x$ the
corresponding primitive idempotent, $S_x$ the corresponding simple
$A$-module, and by $P_x$ and $I_x$ the corresponding indecomposable
projective and injective A-module, respectively.

Let $A$ be an algebra. A module $T_A$ is called a \textit{tilting
module} \cite{hr} if $pd\: T_A \leq 1$; Ext$^1_A(T; T) = 0$ and
the number of isomorphism classes of indecomposable summands of
$T$ equals the rank of the Grothendieck group $K_0(A)$. An algebra
$A$ is called \textit{tilted of type $Q$} if it is the
endomorphism algebra of a tilting $kQ$-module.

An algebra $A$ is called \textit{quasitilted} if $\mathrm{gl.dim.}\:
A \leq 2$ and, for each indecomposable module $M_A$, we have $pd\:M
\leq 1$ or $id\:M \leq 1$ (see \cite{hrs}). It follows from
\cite{hrs}, that tilted algebras are a subclass of quasitilted
algebras, and if a representation finite algebra is quasitilted,
then it is tilted.


\section{On the relations of a quasitilted
algebra}\label{section qt}

The main objective of this section is  the interaction between the
relations of a quasitilted algebra. More precisely, for this class
of algebras we are going to study how the relations can interact.

We begin by introducing the concept of bound consecutive relations
inspired by the ones introduced in \cite{ar}. Observe that the pair
of relations defined in \cite{ar} is a subclass of the bound
consecutive relations defined in the following.

\begin{defn}
The relations $\rho_1$  and $\rho_2$  are called \textit{bound
consecutive relations} if  there is a walk $\varpi$ between
$t(\rho_1)$ and $s(\rho_2)$ such that $\varpi$ does not contain
zero subpaths.
\end{defn}

For example, if $A=kQ/I$ is the algebra given by the quiver

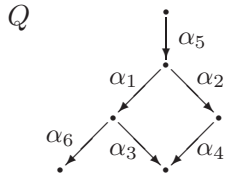
\begin{picture}(0,70)

      \put(160,23){
         \put(-40,35){$Q$}
      \put(0,0){\circle*{2}}
      \put(20,20){\circle*{2}}
      \put(20,-20){\circle*{2}}
      \put(-20,-20){\circle*{2}}
      \put(40,0){\circle*{2}}
      \put(20,40){\circle*{2}}
      \put(18,18){\vector(-1,-1){16}}
      \put(2,-2){\vector(1,-1){16}}
      \put(22,18){\vector(1,-1){16}}
      \put(38,-2){\vector(-1,-1){16}}
      \put(20,38){\vector(0,-1){16}}
      \put(-2,-2){\vector(-1,-1){16}}
          \put(4,15){\HVCenter{\small $\alpha_1$}}
          \put(37,15){\HVCenter{\small $\alpha_2$}}
          \put(4,-12){\HVCenter{\small $\alpha_3$}}
          \put(37,-12){\HVCenter{\small $\alpha_4$}}
          \put(30,30){\HVCenter{\small $\alpha_5$}}
          \put(-20,-7){\HVCenter{\small $\alpha_6$}}
    }
\end{picture}

\noindent bound by $\alpha_5 \alpha_1 \alpha_3 = 0$ and $\alpha_5
\alpha_1 \alpha_6 = 0$. These relations are bound consecutive for
the walk $\varpi= \alpha_4 \alpha_2 \alpha_5$:

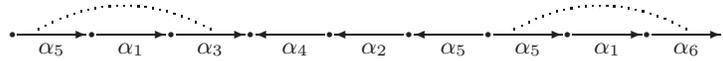
\begin{picture}(0,25)

      \put(80,-10){
      \put(-30,20){\circle*{2}}
      \put(0,20){\circle*{2}}
      \put(30,20){\circle*{2}}
      \put(60,20){\circle*{2}}
      \put(90,20){\circle*{2}}
      \put(120,20){\circle*{2}}
      \put(150,20){\circle*{2}}
      \put(180,20){\circle*{2}}
      \put(210,20){\circle*{2}}
      \put(240,20){\circle*{2}}
      \put(-28,20){\vector(1,0){26}}
          \put(-15,15){\HVCenter{\footnotesize $\alpha_5$}}
      \put(2,20){\vector(1,0){26}}
          \put(15,15){\HVCenter{\footnotesize $\alpha_1$}}
      \put(32,20){\vector(1,0){26}}
          \put(45,15){\HVCenter{\footnotesize $\alpha_3$}}
      \put(88,20){\vector(-1,0){26}}
          \put(77,15){\HVCenter{\footnotesize $\alpha_4$}}
      \put(118,20){\vector(-1,0){26}}
          \put(107,15){\HVCenter{\footnotesize $\alpha_2$}}
      \put(148,20){\vector(-1,0){26}}
          \put(137,15){\HVCenter{\footnotesize $\alpha_5$}}
      \put(152,20){\vector(1,0){26}}
          \put(165,15){\HVCenter{\footnotesize $\alpha_5$}}
      \put(182,20){\vector(1,0){26}}
          \put(195,15){\HVCenter{\footnotesize $\alpha_1$}}
      \put(212,20){\vector(1,0){26}}
          \put(225,15){\HVCenter{\footnotesize $\alpha_6$}}
            \qbezier[24](-20,22)(12,40)(45,22)
            \qbezier[24](160,22)(192,40)(225,22)
    }
\end{picture}\vspace{0.1in}

The main theorem of this section, states that any algebra of
global dimension two that does not contain bound consecutive
relations is quasitilted.

\begin{thm}\label{thm:vuelta}
Let $A$ be an algebra with $\mathrm{gl.dim.}$$\: A \leq 2$.
If there are no bound consecutive relations in $A$, then $A$ is
quasitilted.
\end{thm}

\begin{proof} Suppose $A$ is not a quasitilted algebra. Then, there is an indecomposable $A$-module $M$
such that $pd\: M = id\: M = 2$.

Let

  \begin{picture}(50,30)
      \put(65,-10){
      \put(0,20){\HVCenter{\small $0$}}
      \put(40,20){\HVCenter{\small $P_2$}}
      \put(80,20){\HVCenter{\small $P_1$}}
      \put(120,20){\HVCenter{\small $P_0$}}
      \put(160,20){\HVCenter{\small $M$}}
      \put(200,20){\HVCenter{\small $0$}}
      \put(5,20){\vector(1,0){26}}
      \put(48,20){\vector(1,0){24}}
      \put(88,20){\vector(1,0){24}}
      \put(128,20){\vector(1,0){24}}
      \put(168,20){\vector(1,0){24}}
      \put(60,28){\HVCenter{\small $f_2$}}
      \put(100,28){\HVCenter{\small $f_1$}}
      \put(140,28){\HVCenter{\small $f_0$}}
    }
  \end{picture}

\noindent and

  \begin{picture}(50,30)
      \put(65,-5){
      \put(0,20){\HVCenter{\small $0$}}
      \put(40,20){\HVCenter{\small $M$}}
      \put(80,20){\HVCenter{\small $I_0$}}
      \put(120,20){\HVCenter{\small $I_1$}}
      \put(160,20){\HVCenter{\small $I_2$}}
      \put(200,20){\HVCenter{\small $0$}}
      \put(5,20){\vector(1,0){26}}
      \put(48,20){\vector(1,0){24}}
      \put(88,20){\vector(1,0){24}}
      \put(128,20){\vector(1,0){24}}
      \put(168,20){\vector(1,0){24}}
    }
  \end{picture}

\noindent be the minimal projective and injective resolutions of
$M$, respectively.

Any indecomposable sum $P_2(S_b)$ of $P_2$ gives rise to a diagram
$\begin{array}{ccccc} P_2(S_b) & \longrightarrow & P_1 &
\longrightarrow & P_0  \  \end{array}$ with zero composition. Then
there is a sum $P_0 (S_j)$ of $P_0$ and a relation
$\begin{array}{ccccc}   P_2(S_b) & \longrightarrow & P_1 &
\longrightarrow & P_0(S_j).  \
  \end{array}$ Consequently, there is a relation $\rho_2: j \rightsquigarrow b$
starting at the vertex $j$, where $S_j \in \mbox{Top}\: M$.

Dually, using the minimal injective resolution of $M$, there is $S_i
\in \mbox{Soc}\: M$ and a relation $\rho_1$ ending at the vertex
$i$,

Let $S_k \in \mbox{Top}\: M$ such that $S_i \in f_0(P_0(S_k))$. Then
there exists a nonzero morphism $h: P_0(S_i) \longrightarrow
f_0(P_0(S_k))$ and  an epimorphism $g: P_0(S_k) \longrightarrow
f_0(P_0(S_k))$, i.e. the situation is as follows:

  \begin{picture}(50,65)
      \put(65,10){
      \put(160,-5){\HVCenter{\small $P_0(S_k)$}}
      \put(80,40){\HVCenter{\small $P_0(S_i)$}}
      \put(160,40){\HVCenter{\small $f_0(P_0(S_k))$}}
      \put(100,40){\vector(1,0){30}}
      \put(115,48){\HVCenter{\small $h$}}
      \put(160,5){\vector(0,1){24}}
      \put(170,17){\HVCenter{\small $g$}}
            \qbezier[20](80,30)(110,15)(140,0)
      \put(142,-1){\vector(2,-1){0.1}}
      \put(110,6){\HVCenter{\small $f$}}
    }
  \end{picture}

Then there is a nonzero morphism $f: P_0(S_i) \longrightarrow
P_0(S_k)$ such that $gf = h$. Therefore, there is a nonzero path
$\gamma: k \rightsquigarrow i$.

Since $M$ is an indecomposable module, then $\mbox{sop}\: M$ is
connected. We consider $Q'_0 = (\mbox{sop}\: M)_0$ and $Q'_1 = \{
\alpha \in Q_1 \mbox{ such that } M(\alpha) \neq 0 \}$.
Consequently, $Q' = (Q'_0, Q'_1)$ is a connected quiver and the
vertices associated with the simple $A$-modules of $\mbox{Top}\: M$
are sources of $Q'$. Given two sources $j$ and $k$, there is a walk
of nonzero paths that unites them, (\cite{ars}, pp. 45). Therefore,
$\rho_1$ y $\rho_2$ are  bound consecutive relations
\end{proof}

We want to remark that the algebra $A$, given in the previous
example, is not quasitilted.

We illustrated with the following example that the reciprocal of
Theorem \ref{thm:vuelta}  is generally not true. The algebra $B$
given by the following quiver with relations:

  \begin{picture}(50,80)
      \put(150,15){
      \put(-10,20){\HVCenter{\small $1$}}
      \put(20,20){\HVCenter{\small $2$}}
      \put(20,50){\HVCenter{\small $4$}}
      \put(50,20){\HVCenter{\small $5$}}
      \put(20,-10){\HVCenter{\small $3$}}
      \put(80,20){\HVCenter{\small $6$}}
      \put(15,20){\vector(-1,0){20}}
      \put(24,20){\vector(1,0){20}}
      \put(54,20){\vector(1,0){20}}
      \put(20,44){\vector(0,-1){18}}
      \put(20,-4){\vector(0,1){18}}
            \qbezier[20](22,42)(22,22)(70,22)
            \qbezier[15](18,0)(18,18)(0,18)
    }
  \end{picture}

is a tilted algebra and the given relations are bound consecutive relations, where $\varpi:\ 1 \longleftarrow 2
\longleftarrow 4$ is the path.

Next, we give necessary conditions for an algebra of global
dimension two to be quasitilted. To this end, we introduce the
following notations.

Let $A=kQ/I$ be an algebra and $\rho= \sum^m_{i=1} \lambda_i
\alpha_{i_1}\alpha_{i_2} \ldots \alpha_{i_{r_i}}$ be a relation of
$A$, with $\alpha_{i_j} \in Q_1$. We will note with
$\mathbf{Q}_{\rho}$ the subquiver of $Q$ induced by the arrows of
$\rho$, i.e.  $(\mathbf{Q}_{\rho})_1 = \{ \alpha_{i_j}:
1~\leq~i~\leq~m , 1 \leq j \leq r_i\}$ \  and \
 $(\mathbf{Q}_{\rho})_0 = \{ s(\alpha_{i_j}): \alpha_{i_j} \in
(\mathbf{Q}_{\rho})_1 \} \cup \{ t(\alpha_{i_j}):
\alpha_{i_j}~\in~(\mathbf{Q}_{\rho})_1 \}$.

If $A=kQ/I$ is a schurian strongly simply connected algebra and
 $\rho_1, \rho_2$ are bound consecutive relations for the walk
$$\varpi = \gamma_1\gamma_2\cdots\gamma_r: t(\rho_1)=i
\rightsquigarrow s(\rho_2)=  j$$ with nonzero maximal paths $
\gamma_h $ and $\mathrm{long}$$(\gamma_h)\neq 0$, for $h= 2, \ldots
, r-1$., then we note with $\overline{\mbox{sop}} \: \varpi$ the
support of the union of all the parallel paths of $\gamma_h$ (which
are all equal and nonzero in $A$).

The following result has a strong relation with  Lemma 2.1, given by
\mbox{I. Assem} and M. J. Redondo in \cite{ar}. They proved that a
schurian tilted  algebra $A$ does not contain a certain class of
bound consecutive relations. In the case of schurian strongly simply
connected quasitilted algebras, we show that these algebras do not
have the bounded consecutive relations defined in this work.

\begin{thm}\label{thm:ct1}
Let $A= kQ/I$ be a schurian strongly simply connected algebra with
$\mathrm{gl.dim.}$$\: A = 2$. If $A$ is quasitilted, then there are
no bound consecutive relations $\rho_1$ y $\rho_2$ with walk
$\varpi$ such that $\mathrm{\overline{sop}}$$\; \varpi \; \cap \;
\mathbf{Q}$$_{\rho_1} = \{t(\rho_1)=i\}$ \  and \
$\mathrm{\overline{sop}}$$\; \varpi \; \cap \; \mathbf{Q}$$_{\rho_2}
= \{s(\rho_2)=j\}$.
\end{thm}

\begin{proof} Since $A$ is schurian strongly simply connected,  we can define a representation
$M$ whose support is  $\overline{\mbox{sop}}\; \varpi$ and such that
$M(x)= k$, for all $ x \in (\mbox{sop}\: M )_0$ and $M(\varepsilon)=
Id_k$, for all $\varepsilon \in (\mbox{sop}\: M )_1$. We will show
that the projective dimension of $M$ is 2.

Consider $$P_1(M) \stackrel{f_1}{\longrightarrow} P_0(M)
\stackrel{f_0}{\longrightarrow} M \longrightarrow 0$$ The module $M$
is not projective because $P_j \in \mbox{add} \: P_0(M)$ and
$\mbox{sop} \: P_j \nsubseteq \mbox{sop} \: M$.

If $\rho_2 = \beta \delta_2$, where $\beta \in Q_1$ and $\delta_2$
is a nonzero path in $A$, then $f_0(\beta) = 0$ because
$\mathrm{\overline{sop}}$$\; \varpi \; \cap \; \mathbf{Q}$$_{\rho_2}
= \{s(\rho_2)=j\}$. Then, $P_{t(\beta)}$ is a direct summand of
$P_1(M)$ because the composition $P_{t(\beta)}
\stackrel{\beta\cdot}{\longrightarrow} P_j
\stackrel{f_0}{\longrightarrow} M$  is zero. Therefore, $Im\:
\beta\cdot \subseteq \mbox{Ker}\: f_0$.

As $\beta \notin \mbox{rad Ker}\: f_0$, we have that  $Im\:
\beta\cdot \nsubseteq \mbox{rad Ker}\: f_0$. This proves that
$P_{t(\beta)} \stackrel{\beta\cdot}{\longrightarrow} P_j$ is a sum
of the projective presentation of $M$.

Finally, as $0 = \rho_2 = \beta \delta_2 = \beta\cdot (\delta_2)$,
it follows that $0 \neq \delta_2 \in \mbox{Ker}\; \beta \cdot$.
Consequently, $P_2 (M)\neq 0$.

If $\rho_2 = \beta \delta_2 + \beta' \delta'_2$, where $\beta,
\beta' \in Q_1$ and $\delta_2, \delta'_2$ are  nonzero paths in $A$.
Then, analogously to the previous case proof that $P_{t(\beta)}
\oplus P_{t(\beta')} \stackrel{(\beta\cdot, \beta'\cdot
)}{\longrightarrow} P_j$ is a sum of projective presentation of
$M$. Then, $0 = \rho_2 = \beta \delta_2 + \beta'\delta'_2 = (\beta\cdot,\beta'\cdot) \tiny{\left(\begin{array}{c} \delta_2 \\
\delta'_2 \end{array}\right)}$. Therefore, $P_2 (M)\neq 0$.

By duality, considering the minimal injective co-resolution of $M$,
it follows that  $id\: M = 2$. \end{proof}

\begin{exmp}
Let $A= kQ/I$ be the algebra given by:

  
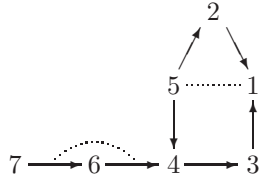
\begin{picture}(50,80)
      \put(150,-10){
      \put(-10,20){\HVCenter{\small $7$}}
      \put(20,20){\HVCenter{\small $6$}}
      \put(50,20){\HVCenter{\small $4$}}
      \put(80,20){\HVCenter{\small $3$}}
      \put(50,50){\HVCenter{\small $5$}}
      \put(80,50){\HVCenter{\small $1$}}
      \put(65,78){\HVCenter{\small $2$}}
      \put(-5,20){\vector(1,0){20}}
      \put(24,20){\vector(1,0){20}}
      \put(54,20){\vector(1,0){20}}
      \put(50,44){\vector(0,-1){18}}
      \put(80,26){\vector(0,1){18}}
      \put(52,56){\vector(1,2){8}}
      \put(70,72){\vector(1,-2){8}}
            \qbezier[10](55,50)(65,50)(75,50)
            \qbezier[15](5,22)(20,35)(35,22)
    }
  \end{picture}

The two relations of this algebra are minimal bound consecutive
relations  with  a walk $\varpi: \  4 \longleftarrow 5$:

As the $A$-module $M$, whose representation is:

  
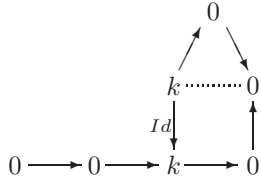
\begin{picture}(50,80)
      \put(150,-10){
      \put(-10,20){\HVCenter{\small $0$}}
      \put(20,20){\HVCenter{\small $0$}}
      \put(50,20){\HVCenter{\small $k$}}
      \put(80,20){\HVCenter{\small $0$}}
      \put(50,50){\HVCenter{\small $k$}}
      \put(80,50){\HVCenter{\small $0$}}
      \put(65,78){\HVCenter{\small $0$}}
      \put(45,35){\HVCenter{\tiny $Id$}}
      \put(-5,20){\vector(1,0){20}}
      \put(24,20){\vector(1,0){20}}
      \put(54,20){\vector(1,0){20}}
      \put(50,44){\vector(0,-1){18}}
      \put(80,26){\vector(0,1){18}}
      \put(52,56){\vector(1,2){8}}
      \put(70,72){\vector(1,-2){8}}
            \qbezier[10](55,50)(65,50)(75,50)
    }
  \end{picture}

\noindent has projective dimension and injective dimension equal to
2. Then, by the previous theorem, $A$ is not a quasitilted algebra.
\end{exmp}

The next question is to study what happens when a quasitilted
algebra has a bound consecutive relation by a walk  $\varpi$ such
that $\overline{\mbox{sop}}\; \varpi$ intersects the quiver of one
of the relations in at least one arrow. In the next theorem we
answer the question posed above, to schurian strongly simply
connected algebras.

\begin{thm}\label{thm:ct2}
Let $A= kQ/I$ be a schurian strongly simply connected algebra and
quasitilted such that there are  $\rho_1$ and $\rho_2$ bound
consecutive relations by the walk $\varpi$. Then, either
$\mathrm{\overline{sop}}\: \varpi \cap Q_{\rho_1}$ \  or \
$\mathrm{\overline{sop}}\: \varpi \cap Q_{\rho_2}$ is a maximal
subpath of $\rho_1$ or $\rho_2$, respectively.
\end{thm}

\begin{proof} Suppose that for $l= 1,2$, \
$\mathrm{\overline{sop}}\: \varpi$ \  and \  $Q_{\rho_l}$ do not
intersect in a maximal subpath of  $\rho_l$. Let $M$  be the
representation whose support is $\overline{\mbox{sop}}\; \varpi$ and
such that $M(x)= k$, for all $x \in (\mbox{sop}\: M )_0$ and
$M(\varepsilon)= Id_k$, for all $\varepsilon \in (\mbox{sop}\: M
)_1$. If we prove that the projective dimension of $M$ is 2, then by
duality $di\: M = 2$.

As $M$ is not a projective module, we can consider $$0 \neq P_1(M)
\stackrel{f_1}{\longrightarrow} P_0(M)
\stackrel{f_0}{\longrightarrow} M \longrightarrow 0.$$

For $\rho_2 =  \delta_2 \beta$, where $\beta \in Q_1$ and $\delta_2$
is a nonzero path in $A$, must be $f_0(\delta_2) = 0$. Writing
$\delta_2 = \psi_2 \varphi_2$, with $\psi_2$ minimal path in
$\delta_2$ such that $f_0(\psi_2)= 0$, results that the composition
$P_{t(\psi_2)} \stackrel{\psi_2\cdot}{\longrightarrow} P_j
\stackrel{f_0}{\longrightarrow} M$ is zero. Therefore, $Im\:
\psi_2\cdot \subseteq \mbox{Ker}\: f_0$.

Suppose that $\psi_2 \in \mbox{rad Ker}\: f_0$. Then there exist
non-trivial paths $\lambda_l$  such that $\psi_2 = \sum \mu_l
\lambda_l $, with $\mu_l \in \mbox{Ker}\: f_0$. As $A$ schurian
strongly simply connected, $\psi_2 = a \lambda_1 \mu_1$, with $a$
scalar. But $f_0(\mu_1)=0$ and $l(\mu_1)< l(\psi_2)$, which
contradicts the minimality of  $\psi_2$. Then, $Im\: \psi_2\cdot
\nsubseteq \mbox{rad Ker}\: f_0$. Consequently, $P_{t(\psi_2)}
\stackrel{\psi_2\cdot}{\longrightarrow} P_j$  is a sum of projective
presentation of $M$.

Since $0 = \rho_2 = \psi_2 \varphi_2 \beta = \psi_2\cdot (\varphi_2
\beta)$, Then $0 \neq \varphi_2\beta \in \mbox{Ker}\; \psi_2 \cdot
$. Therefore,  $P_2 (M)\neq 0$.

For $\rho_2 =  \delta_2 + \delta'_2$, where $\delta_2, \delta'_2$
 are nonzero paths in $A$, it follows that $f_0(\delta_2) = f_0(\delta'_2) =
0$. If we write $\delta_h = \psi_h \varphi_h$, with $\psi_h$ minimal
path in $\delta_h$ such that $f_0(\psi_h)= 0$, for $h = 1,2$, is
proved similarly to the previous case that  $P_{t(\psi_2)} \oplus
P_{t(\psi'_2)} \stackrel{(\psi_2\cdot, \psi'_2\cdot
)}{\longrightarrow} P_j$ is a sum of projective presentation of $M$.

Finally, since $0 = \rho_2 = \delta_2 + \delta'_2 = \psi_2 \varphi_2
+ \psi'_2 \varphi'_2=
(\psi_2\cdot,\psi'_2\cdot) \tiny{\left(\begin{array}{c} \varphi_2 \\
\varphi'_2 \end{array}\right)}$ we have that $P_2 (M)\neq 0$.
\end{proof}

\begin{exmp}
Let $A=kQ/I$ be the tilted algebra given by the following quiver
with relations:

  
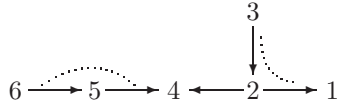
\begin{picture}(50,60)
      \put(130,-10){
      \put(-10,20){\HVCenter{\small $6$}}
      \put(20,20){\HVCenter{\small $5$}}
      \put(50,20){\HVCenter{\small $4$}}
      \put(80,20){\HVCenter{\small $2$}}
      \put(110,20){\HVCenter{\small $1$}}
      \put(80,50){\HVCenter{\small $3$}}
      \put(-5,20){\vector(1,0){20}}
      \put(24,20){\vector(1,0){20}}
      \put(76,20){\vector(-1,0){20}}
      \put(80,44){\vector(0,-1){18}}
      \put(84,20){\vector(1,0){20}}
            \qbezier[10](83,40)(83,25)(95,23)
            \qbezier[15](0,22)(20,35)(35,22)
    }
  \end{picture}

Note that the two monomial relations having $A$ are bound
consecutive by the path $\varpi: 4 \longleftarrow 2 \longleftarrow
3$, that is:

  \begin{picture}(10,30)
      \put(90,-15){
      \put(-10,20){\HVCenter{\small $6$}}
      \put(20,20){\HVCenter{\small $5$}}
      \put(50,20){\HVCenter{\small $4$}}
      \put(80,20){\HVCenter{\small $2$}}
      \put(110,20){\HVCenter{\small $3$}}
      \put(140,20){\HVCenter{\small $2$}}
      \put(170,20){\HVCenter{\small $1$}}
      \put(-6,20){\vector(1,0){20}}
      \put(24,20){\vector(1,0){20}}
      \put(76,20){\vector(-1,0){20}}
      \put(106,20){\vector(-1,0){20}}
      \put(114,20){\vector(1,0){20}}
      \put(144,20){\vector(1,0){20}}
            \qbezier[10](155,22)(140,35)(120,22)
            \qbezier[10](35,22)(20,35)(0,22)
    }
  \end{picture}

\noindent However, $\overline{\mbox{sop}}\: \varpi$ intersects the
relation  $\rho: 3\longrightarrow 2 \longrightarrow 1$ in the path
$3 \longrightarrow 2$, which is maximal in $\rho$.
\end{exmp}

The condition of Theorem \ref{thm:ct2} is not a sufficient condition
to ensure that one is quasitilted algebra. In fact,
let $A$ be the algebra given by the following quiver with relations:

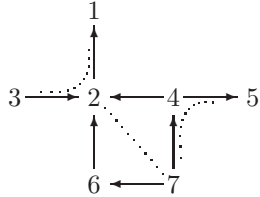
\begin{picture}(0,80)

      \put(160,15){
      \put(-30,20){\HVCenter{\small $3$}}
      \put(0,20){\HVCenter{\small $2$}}
      \put(30,20){\HVCenter{\small $4$}}
      \put(0,53){\HVCenter{\small $1$}}
      \put(60,20){\HVCenter{\small $5$}}
      \put(0,-13){\HVCenter{\small $6$}}
      \put(30,-13){\HVCenter{\small $7$}}
      \put(-26,20){\vector(1,0){20}}
      \put(26,20){\vector(-1,0){20}}
      \put(34,20){\vector(1,0){20}}
      \put(26,-13){\vector(-1,0){20}}
      \put(0,27){\vector(0,1){20}}
      \put(0,-7){\vector(0,1){20}}
      \put(30,-7){\vector(0,1){20}}
            \qbezier[10](26,-9)(15,4)(4,16)
            \qbezier[10](-20,22)(0,20)(-2,38)
            \qbezier[10](33,-3)(30,20)(45,18)
    }
\end{picture}

and we will consider the following monomial relations:
$\rho_1 = 7 \longrightarrow 4 \longrightarrow 5$ \  \  and \   \
$\rho_2 = 3 \longrightarrow 2 \longrightarrow 1$.
Then, $\rho_1$ and $\rho_2$ are bound consecutive relations by the
walk \mbox{$\varpi: 5 \longleftarrow 4 \longrightarrow 2
\longleftarrow 3$}:

  \begin{picture}(10,30)
      \put(85,-15){
      \put(-10,20){\HVCenter{\small $7$}}
      \put(20,20){\HVCenter{\small $4$}}
      \put(50,20){\HVCenter{\small $5$}}
      \put(80,20){\HVCenter{\small $4$}}
      \put(110,20){\HVCenter{\small $2$}}
      \put(140,20){\HVCenter{\small $3$}}
      \put(170,20){\HVCenter{\small $2$}}
      \put(200,20){\HVCenter{\small $1$}}
      \put(-6,20){\vector(1,0){20}}
      \put(24,20){\vector(1,0){20}}
      \put(76,20){\vector(-1,0){20}}
      \put(84,20){\vector(1,0){20}}
      \put(136,20){\vector(-1,0){20}}
      \put(144,20){\vector(1,0){20}}
      \put(174,20){\vector(1,0){20}}
            \qbezier[10](185,22)(170,35)(150,22)
            \qbezier[10](35,22)(20,35)(0,22)
    }
  \end{picture}

It is clear that $\mathrm{\overline{sop}}\: \varpi$ intersects
$Q_{\rho_1}$ and $Q_{\rho_2}$ in maximal subpaths of  $\rho_1$ and
$\rho_2$, respectively.

However, the $A$-module $ M $ given by the following representation

\begin{picture}(0,80)

      \put(160,15){
      \put(-30,20){\HVCenter{\small $k$}}
      \put(0,20){\HVCenter{\small $k$}}
      \put(30,20){\HVCenter{\small $k$}}
      \put(0,53){\HVCenter{\small $0$}}
      \put(60,20){\HVCenter{\small $k$}}
      \put(0,-13){\HVCenter{\small $0$}}
      \put(30,-13){\HVCenter{\small $0$}}
      \put(-26,20){\vector(1,0){20}}
      \put(25,20){\vector(-1,0){20}}
      \put(34,20){\vector(1,0){20}}
      \put(25,-13){\vector(-1,0){20}}
      \put(0,27){\vector(0,1){20}}
      \put(0,-7){\vector(0,1){20}}
      \put(30,-7){\vector(0,1){20}}
      \put(16,-18){\HVCenter{\tiny $0$}}
      \put(-5,2){\HVCenter{\tiny $0$}}
      \put(35,2){\HVCenter{\tiny $0$}}
      \put(-5,36){\HVCenter{\tiny $0$}}
      \put(-17,25){\HVCenter{\tiny $Id$}}
      \put(16,25){\HVCenter{\tiny $Id$}}
      \put(43,25){\HVCenter{\tiny $Id$}}
    }
\end{picture}\vspace{.05in}

\noindent is indecomposable and  \   $pd \: M = id \: M = 2$.
Therefore, $A$ is not quasitilted algebra.

\section{On the quiver with relations of a quasitilted algebra}

We study the quiver with relations of a quasitilted algebra. We start recalling the definition of a critical algebra given in \cite{bft}.

\begin{defn}\label{def} Let  $B$ be an  algebra. We say that  $B$ is
\textit{critical} if it satisfies the following properties:

\begin{enumerate}
    \item[i)] $B$ has a unique source $i$ and a unique sink  $j$,
    \item[ii)] $pd\: S_i = id\: S_j = 3$ and if $S$ is a different simple $B$-module we have that $pd\:  S \leq 2 $ and $id\: S \leq
    2$,
    \item[iii)] Let $0 \rightarrow \textsf{Q}_3 \rightarrow \textsf{Q}_2
    \rightarrow \textsf{Q}_1 \rightarrow \textsf{Q}_0 \rightarrow S_i
    \rightarrow 0$ be the minimal projective resolution  of the simple $S_i$ and let $\textsf{Q} = \oplus^{3}_{k=0} \textsf{Q}_k$.
    Then all  indecomposable projective $B$-modules are in $\mathrm{add}$$\: \textsf{Q}$, and each  indecomposable projective $B$-module
    is a direct summand of exactly one $\textsf{Q}_k$,
    \item[iv)] Let $0 \rightarrow S_j \rightarrow I_0 \rightarrow I_1 \rightarrow I_2 \rightarrow I_3
    \rightarrow 0$ be the minimal co-injective resolution  of simple  $S_j$ and let $I = \oplus^{3}_{k=0} I_k$. Then all
    indecomposable injective $B$-modules are in $\mathrm{add}$$\: I$, and each indecomposable injective $B$-module
    is a direct summand of exactly one $I_k$.

    \item[v)] $B$ does not contain any proper full subcategory that verifies i), ii), iii) and iv).
\end{enumerate}
\end{defn}

We have the following corollary.

\begin{coro} Let A be a strongly simply connected schurian algebra. If the following
conditions hold:
\begin{enumerate}
\item[a)] A does not contain a proper full subcategory B such that B is a critical algebra.
\item[b)] A has no bounded consecutive relations.
\end{enumerate}
Then A is a quasitilted algebra.
\end{coro}

\begin{proof}
It follows from a) and \cite{bft} that $\mathrm{gl.dim.} A \leq 2$.
Appling Theorem \ref{thm:vuelta} we get that $A$ is quasitilted
algebra.
\end{proof}

\section{Applications}

\subsection{Tilted algebras of Dynkin type}

In this section we provide a method to obtain a complete classification of tilted
algebras of Dynkin type. In particular we give the classification
of tilted algebras of type $E_6$.

The combinatorial description of how the iterated tilted algebra
of Dynkin type can be obtained from a trivial extension of finite
representation type was studied in \cite{f}.  One of the essential
tools in this description is the notion of admissible cuts of
trivial extensions (see \cite{f} and \cite{fp}).

The next step is to select the iterated tilted algebras with
global dimension at most two, using \cite{bft}. For
this propose we use the Theorem \ref{thm:vuelta}, Theorem
\ref{thm:ct1} and Theorem \ref{thm:ct2}.

In this way, we could obtain a classification of the tilted algebras of
Dynkin type. We now give the classification of Dynkin
type $E_6$.

In the following, we will consider figures obtained from a quiver
considering some arrows without orientation. Moreover, it is said
that a quiver is associated with a figure if it is obtained giving
an orientation to these arcs.

\begin{thm}\label{thm:te6}
Let $B$ be an algebra. Then $B$ is tilted of Dynkin type $E_6$ if
and only if
\begin{enumerate}
    \item[i.] $B$ is hereditary, or
    \item[ii.] $B$ or $B^{op}$ is isomorphic to one of the following:
\end{enumerate}
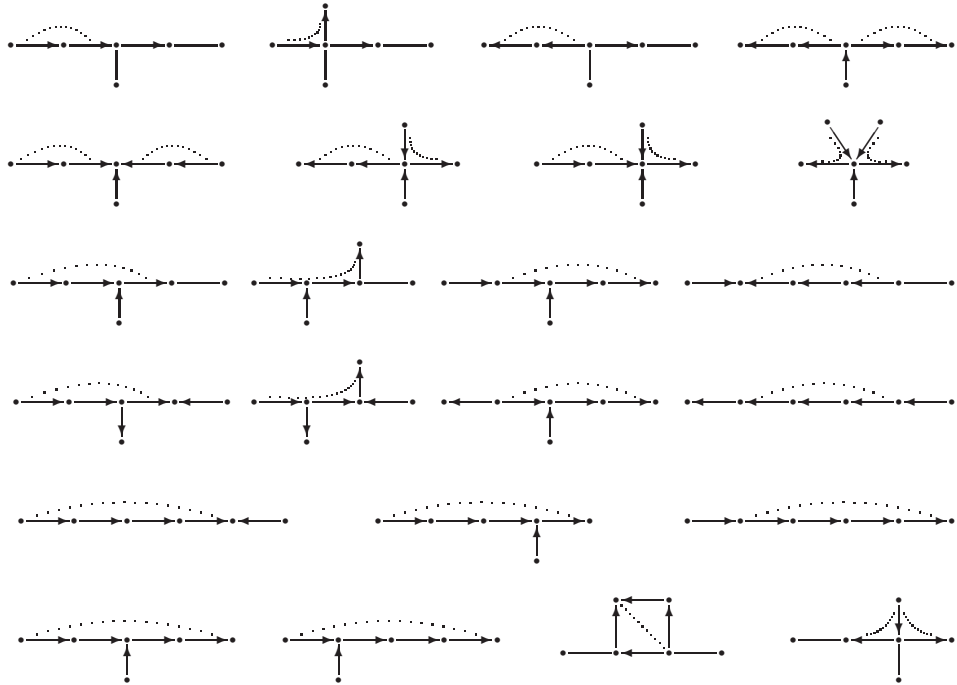
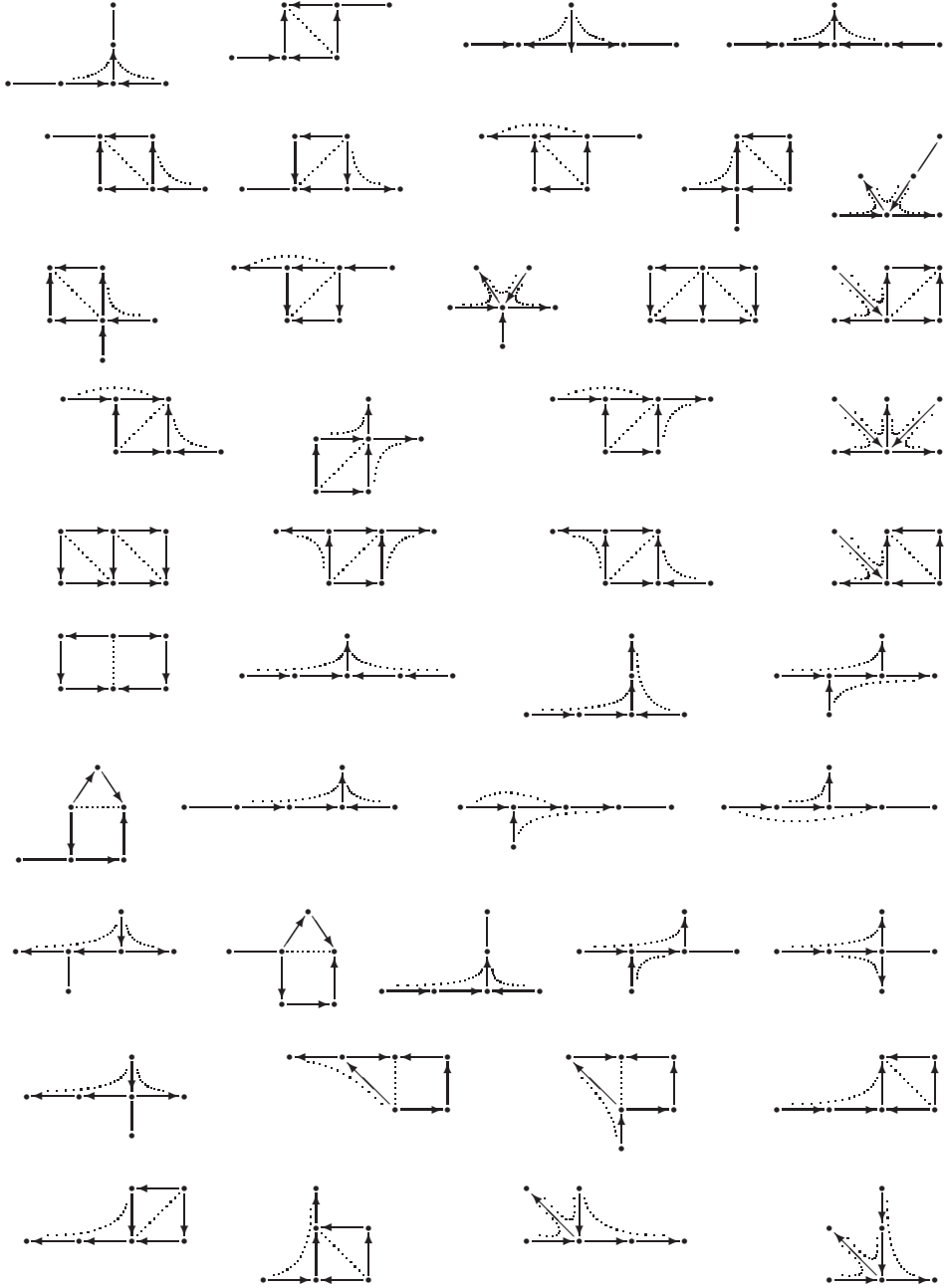
\vspace{.05in}

  \begin{picture}(50,30)

      \put(21,-5){
      \put(0,20){\circle*{2}}
      \put(-20,20){\circle*{2}}
      \put(20,20){\circle*{2}}
      \put(60,20){\circle*{2}}
      \put(40,20){\circle*{2}}
      \put(20,5){\circle*{2}}
      \put(-18,20){\vector(1,0){16}}
      \multiput(2,20)(20,0){2}{\vector(1,0){16}}
      \put(20,7){\line(0,2){11}}
      \put(42,20){\line(1,0){16}}
             \qbezier[12](-14,22)(0,32)(10,22)
      }

      \put(100,-5){
      \put(0,20){\circle*{2}}
      \put(20,35){\circle*{2}}
      \put(20,20){\circle*{2}}
      \put(60,20){\circle*{2}}
      \put(40,20){\circle*{2}}
      \put(20,5){\circle*{2}}
      \put(20,22){\vector(0,1){11}}
      \multiput(2,20)(20,0){2}{\vector(1,0){16}}
      \put(20,7){\line(0,2){11}}
      \put(42,20){\line(1,0){16}}
             \qbezier[10](6,22)(18,22)(18,28)
      }

      \put(200,-5){
      \put(0,20){\circle*{2}}
      \put(-20,20){\circle*{2}}
      \put(20,20){\circle*{2}}
      \put(60,20){\circle*{2}}
      \put(40,20){\circle*{2}}
      \put(20,5){\circle*{2}}
      \put(22,20){\vector(1,0){16}}
      \multiput(18,20)(-20,0){2}{\vector(-1,0){16}}
      \put(20,7){\line(0,2){11}}
      \put(42,20){\line(1,0){16}}
             \qbezier[12](-12,22)(0,32)(14,22)
      }


      \put(217,-5){
      \put(80,20){\circle*{2}}
      \put(60,20){\circle*{2}}
      \put(100,20){\circle*{2}}
      \put(140,20){\circle*{2}}
      \put(120,20){\circle*{2}}
      \put(100,5){\circle*{2}}
      \multiput(98,20)(-20,0){2}{\vector(-1,0){16}}
      \multiput(102,20)(20,0){2}{\vector(1,0){16}}
      \put(100,7){\vector(0,2){11}}
            \qbezier[12](94,22)(80,32)(68,22)
            \qbezier[12](132,22)(120,32)(106,22)
    }
\end{picture}

  \begin{picture}(50,40)
      \put(-209,-5){
      \put(230,20){\circle*{2}}
      \put(210,20){\circle*{2}}
      \put(250,20){\circle*{2}}
      \put(290,20){\circle*{2}}
      \put(270,20){\circle*{2}}
      \put(250,5){\circle*{2}}
      \multiput(288,20)(-20,0){2}{\vector(-1,0){16}}
      \multiput(212,20)(20,0){2}{\vector(1,0){16}}
      \put(250,7){\vector(0,2){11}}
            \qbezier[12](240,22)(230,32)(214,22)
            \qbezier[12](284,22)(270,32)(260,22)
    }

    \put(50,95){
      \put(80,-80){\circle*{2}}
      \put(60,-80){\circle*{2}}
      \put(100,-80){\circle*{2}}
      \put(100,-65){\circle*{2}}
      \put(120,-80){\circle*{2}}
      \put(100,-95){\circle*{2}}
      \multiput(98,-80)(-20,0){2}{\vector(-1,0){16}}
      \put(100,-67){\vector(0,-1){11}}
      \put(102,-80){\vector(1,0){16}}
      \put(100,-93){\vector(0,2){11}}
            \qbezier[12](94,-78)(80,-68)(68,-78)
            \qbezier[10](102,-70)(102,-78)(112,-78)
    }

      \put(-10,95){
      \put(230,-80){\circle*{2}}
      \put(210,-80){\circle*{2}}
      \put(250,-80){\circle*{2}}
      \put(250,-65){\circle*{2}}
      \put(270,-80){\circle*{2}}
      \put(250,-95){\circle*{2}}
      \multiput(212,-80)(20,0){2}{\vector(1,0){16}}
      \put(250,-67){\vector(0,-1){11}}
      \put(252,-80){\vector(1,0){16}}
      \put(250,-93){\vector(0,2){11}}
            \qbezier[12](242,-78)(230,-68)(216,-78)
            \qbezier[10](252,-70)(252,-78)(262,-78)
    }

       \put(250,95){
      \put(50,-80){\circle*{2}}
      \put(60,-64){\circle*{2}}
      \put(70,-80){\circle*{2}}
      \put(80,-64){\circle*{2}}
      \put(90,-80){\circle*{2}}
      \put(70,-95){\circle*{2}}
      \put(61,-66){\vector(2,-3){8}}
      \put(68,-80){\vector(-1,0){16}}
      \put(79,-66){\vector(-2,-3){8}}
      \put(72,-80){\vector(1,0){16}}
      \put(70,-93){\vector(0,2){11}}
            \qbezier[9](79,-70)(70,-80)(82,-79)
            \qbezier[9](61,-70)(70,-80)(58,-79)
    }
\end{picture}

  \begin{picture}(50,40)
%

      \put(22,-5){
      \put(0,20){\circle*{2}}
      \put(-20,20){\circle*{2}}
      \put(20,20){\circle*{2}}
      \put(60,20){\circle*{2}}
      \put(40,20){\circle*{2}}
      \put(20,5){\circle*{2}}
      \put(-18,20){\vector(1,0){16}}
      \multiput(2,20)(20,0){2}{\vector(1,0){16}}
      \put(20,7){\vector(0,2){11}}
      \put(42,20){\line(1,0){16}}
             \qbezier[12](-14,22)(15,32)(30,22)
      }

      \put(93,-5){
      \put(0,20){\circle*{2}}
      \put(40,35){\circle*{2}}
      \put(20,20){\circle*{2}}
      \put(60,20){\circle*{2}}
      \put(40,20){\circle*{2}}
      \put(20,5){\circle*{2}}
      \put(40,22){\vector(0,1){11}}
      \multiput(2,20)(20,0){2}{\vector(1,0){16}}
      \put(20,7){\vector(0,2){11}}
      \put(42,20){\line(1,0){16}}
             \qbezier[15](6,22)(38,20)(38,28)
      }

      \put(185,-5){
      \put(0,20){\circle*{2}}
      \put(-20,20){\circle*{2}}
      \put(20,20){\circle*{2}}
      \put(60,20){\circle*{2}}
      \put(40,20){\circle*{2}}
      \put(20,5){\circle*{2}}
      \put(-18,20){\vector(1,0){16}}
      \multiput(2,20)(20,0){2}{\vector(1,0){16}}
      \put(20,7){\vector(0,2){11}}
      \put(42,20){\vector(1,0){16}}
             \qbezier[12](6,22)(32,32)(52,22)
      }

      \put(277,-5){
      \put(0,20){\circle*{2}}
      \put(-20,20){\circle*{2}}
      \put(20,20){\circle*{2}}
      \put(60,20){\circle*{2}}
      \put(40,20){\circle*{2}}
      \put(80,20){\circle*{2}}
      \put(-18,20){\vector(1,0){16}}
      \multiput(38,20)(-20,0){2}{\vector(-1,0){16}}
      \put(62,20){\line(1,0){16}}
      \put(58,20){\vector(-1,0){16}}
             \qbezier[12](8,22)(32,32)(52,22)
      }
\end{picture}

  \begin{picture}(50,40)
%
      \put(23,-5){
      \put(0,20){\circle*{2}}
      \put(-20,20){\circle*{2}}
      \put(20,20){\circle*{2}}
      \put(60,20){\circle*{2}}
      \put(40,20){\circle*{2}}
      \put(20,5){\circle*{2}}
      \put(-18,20){\vector(1,0){16}}
      \multiput(2,20)(20,0){2}{\vector(1,0){16}}
      \put(20,18){\vector(0,-2){11}}
      \put(58,20){\vector(-1,0){16}}
             \qbezier[12](-14,22)(15,32)(30,22)
      }

      \put(93,-5){
      \put(0,20){\circle*{2}}
      \put(40,35){\circle*{2}}
      \put(20,20){\circle*{2}}
      \put(60,20){\circle*{2}}
      \put(40,20){\circle*{2}}
      \put(20,5){\circle*{2}}
      \put(40,22){\vector(0,1){11}}
      \multiput(2,20)(20,0){2}{\vector(1,0){16}}
      \put(20,18){\vector(0,-2){11}}
      \put(58,20){\vector(-1,0){16}}
             \qbezier[15](6,22)(35,20)(38,28)
      }

      \put(185,-5){
      \put(0,20){\circle*{2}}
      \put(-20,20){\circle*{2}}
      \put(20,20){\circle*{2}}
      \put(60,20){\circle*{2}}
      \put(40,20){\circle*{2}}
      \put(20,5){\circle*{2}}
      \put(-2,20){\vector(-1,0){16}}
      \multiput(2,20)(20,0){2}{\vector(1,0){16}}
      \put(20,7){\vector(0,2){11}}
      \put(42,20){\vector(1,0){16}}
             \qbezier[12](6,22)(32,32)(52,22)
      }

      \put(277,-5){
      \put(0,20){\circle*{2}}
      \put(-20,20){\circle*{2}}
      \put(20,20){\circle*{2}}
      \put(60,20){\circle*{2}}
      \put(40,20){\circle*{2}}
      \put(80,20){\circle*{2}}
      \put(-2,20){\vector(-1,0){16}}
      \multiput(38,20)(-20,0){2}{\vector(-1,0){16}}
      \put(78,20){\vector(-1,0){16}}
      \put(58,20){\vector(-1,0){16}}
             \qbezier[12](8,22)(32,32)(54,22)
      }
\end{picture}

  \begin{picture}(50,40)

      \put(25,-5){
      \put(0,20){\circle*{2}}
      \put(-20,20){\circle*{2}}
      \put(20,20){\circle*{2}}
      \put(60,20){\circle*{2}}
      \put(40,20){\circle*{2}}
      \put(80,20){\circle*{2}}
      \put(-18,20){\vector(1,0){16}}
      \multiput(2,20)(20,0){2}{\vector(1,0){16}}
      \put(78,20){\vector(-1,0){16}}
      \put(42,20){\vector(1,0){16}}
             \qbezier[16](-14,22)(19,32)(52,22)
      }

      \put(160,-5){
      \put(0,20){\circle*{2}}
      \put(-20,20){\circle*{2}}
      \put(20,20){\circle*{2}}
      \put(60,20){\circle*{2}}
      \put(40,20){\circle*{2}}
      \put(40,5){\circle*{2}}
      \put(-18,20){\vector(1,0){16}}
      \multiput(2,20)(20,0){2}{\vector(1,0){16}}
      \put(40,7){\vector(0,2){11}}
      \put(42,20){\vector(1,0){16}}
             \qbezier[16](-14,22)(19,32)(52,22)
      }

      \put(297,-5){
      \put(0,20){\circle*{2}}
      \put(-20,20){\circle*{2}}
      \put(20,20){\circle*{2}}
      \put(60,20){\circle*{2}}
      \put(40,20){\circle*{2}}
      \put(-40,20){\circle*{2}}
      \put(-18,20){\vector(1,0){16}}
      \multiput(2,20)(20,0){2}{\vector(1,0){16}}
      \put(-38,20){\vector(1,0){16}}
      \put(42,20){\vector(1,0){16}}
             \qbezier[16](-14,22)(19,32)(52,22)
      }
\end{picture}

  \begin{picture}(50,40)
      \put(25,-5){
      \put(0,20){\circle*{2}}
      \put(-20,20){\circle*{2}}
      \put(20,20){\circle*{2}}
      \put(60,20){\circle*{2}}
      \put(40,20){\circle*{2}}
      \put(20,5){\circle*{2}}
      \put(-18,20){\vector(1,0){16}}
      \multiput(2,20)(20,0){2}{\vector(1,0){16}}
      \put(20,7){\vector(0,2){11}}
      \put(42,20){\vector(1,0){16}}
             \qbezier[16](-14,22)(19,32)(52,22)
      }

      \put(125,-5){
      \put(0,20){\circle*{2}}
      \put(-20,20){\circle*{2}}
      \put(20,20){\circle*{2}}
      \put(60,20){\circle*{2}}
      \put(40,20){\circle*{2}}
      \put(0,5){\circle*{2}}
      \put(-18,20){\vector(1,0){16}}
      \multiput(2,20)(20,0){2}{\vector(1,0){16}}
      \put(0,7){\vector(0,2){11}}
      \put(42,20){\vector(1,0){16}}
             \qbezier[16](-14,22)(19,32)(52,22)
      }


      \put(210,-10){
      \put(0,20){\circle*{2}}
      \put(20,20){\circle*{2}}
      \put(40,20){\circle*{2}}
      \put(60,20){\circle*{2}}
      \put(40,40){\circle*{2}}
      \put(20,40){\circle*{2}}
      \put(38,20){\vector(-1,0){16}}
      \put(42,20){\line(1,0){16}}
      \put(18,20){\line(-1,0){16}}
      \put(40,22){\vector(0,1){16}}
      \put(38,40){\vector(-1,0){16}}
      \put(20,22){\vector(0,1){16}}
            \qbezier[10](22,38)(30,30)(38,22)
    }

      \put(297,-5){
      \put(0,20){\circle*{2}}
      \put(40,5){\circle*{2}}
      \put(40,20){\circle*{2}}
      \put(60,20){\circle*{2}}
      \put(40,35){\circle*{2}}
      \put(20,20){\circle*{2}}
      \put(38,20){\vector(-1,0){16}}
      \put(40,18){\line(0,-1){11}}
      \put(18,20){\line(-1,0){16}}
      \put(40,33){\vector(0,-1){11}}
      \put(42,20){\vector(1,0){16}}
             \qbezier[10](28,22)(35,23)(38,30)
             \qbezier[10](42,30)(45,23)(52,22)
    }
\end{picture}

  \begin{picture}(50,45)
      \put(3,-15){
      \put(0,20){\circle*{2}}
      \put(20,20){\circle*{2}}
      \put(40,20){\circle*{2}}
      \put(60,20){\circle*{2}}
      \put(40,35){\circle*{2}}
      \put(40,50){\circle*{2}}
      \put(22,20){\vector(1,0){16}}
      \put(58,20){\vector(-1,0){16}}
      \put(18,20){\line(-1,0){16}}
      \put(40,37){\line(0,1){11}}
      \put(40,22){\vector(0,1){11}}
             \qbezier[10](25,22)(35,23)(38,28)
             \qbezier[10](42,28)(45,23)(55,22)
    }


      \put(88,-5){
      \put(0,20){\circle*{2}}
      \put(20,20){\circle*{2}}
      \put(40,20){\circle*{2}}
      \put(60,40){\circle*{2}}
      \put(40,40){\circle*{2}}
      \put(20,40){\circle*{2}}
      \put(38,20){\vector(-1,0){16}}
      \put(42,40){\line(1,0){16}}
      \put(2,20){\vector(1,0){16}}
      \put(40,22){\vector(0,1){16}}
      \put(38,40){\vector(-1,0){16}}
      \put(20,22){\vector(0,1){16}}
            \qbezier[10](22,38)(30,30)(38,22)
    }

      \put(177,-0){
      \put(0,20){\circle*{2}}
      \put(60,20){\circle*{2}}
      \put(40,20){\circle*{2}}
      \put(80,20){\circle*{2}}
      \put(40,35){\circle*{2}}
      \put(20,20){\circle*{2}}
      \put(38,20){\vector(-1,0){16}}
      \put(62,20){\line(1,0){16}}
      \put(2,20){\vector(1,0){16}}
      \put(40,33){\vector(0,-1){16}}
      \put(42,20){\vector(1,0){16}}
             \qbezier[10](28,22)(35,23)(38,30)
             \qbezier[10](42,30)(45,23)(52,22)
    }

      \put(277,-0){
      \put(0,20){\circle*{2}}
      \put(60,20){\circle*{2}}
      \put(40,20){\circle*{2}}
      \put(80,20){\circle*{2}}
      \put(40,35){\circle*{2}}
      \put(20,20){\circle*{2}}
      \put(22,20){\vector(1,0){16}}
      \put(78,20){\vector(-1,0){16}}
      \put(2,20){\vector(1,0){16}}
      \put(40,22){\vector(0,1){11}}
      \put(58,20){\vector(-1,0){16}}
             \qbezier[10](25,22)(35,23)(38,28)
             \qbezier[10](42,28)(45,23)(55,22)
    }
\end{picture}

  \begin{picture}(50,45)

      \put(18,-5){
      \put(0,40){\circle*{2}}
      \put(20,20){\circle*{2}}
      \put(40,20){\circle*{2}}
      \put(60,20){\circle*{2}}
      \put(40,40){\circle*{2}}
      \put(20,40){\circle*{2}}
      \put(38,20){\vector(-1,0){16}}
      \put(18,40){\line(-1,0){16}}
      \put(40,22){\vector(0,1){16}}
      \put(38,40){\vector(-1,0){16}}
      \put(20,22){\vector(0,1){16}}
      \put(58,20){\vector(-1,0){16}}
            \qbezier[10](22,38)(30,30)(38,22)
            \qbezier[10](42,32)(45,23)(54,22)
    }

      \put(92,-5){
      \put(0,20){\circle*{2}}
      \put(20,20){\circle*{2}}
      \put(40,20){\circle*{2}}
      \put(60,20){\circle*{2}}
      \put(40,40){\circle*{2}}
      \put(20,40){\circle*{2}}
      \put(38,20){\vector(-1,0){16}}
      \put(18,20){\line(-1,0){16}}
      \put(40,38){\vector(0,-1){16}}
      \put(38,40){\vector(-1,0){16}}
      \put(20,38){\vector(0,-1){16}}
      \put(42,20){\vector(1,0){16}}
            \qbezier[10](22,22)(30,30)(38,38)
            \qbezier[10](42,34)(43,23)(52,22)
    }

      \put(183,-5){
      \put(0,40){\circle*{2}}
      \put(20,20){\circle*{2}}
      \put(40,20){\circle*{2}}
      \put(60,40){\circle*{2}}
      \put(40,40){\circle*{2}}
      \put(20,40){\circle*{2}}
      \put(38,20){\vector(-1,0){16}}
      \put(18,40){\vector(-1,0){16}}
      \put(40,22){\vector(0,1){16}}
      \put(38,40){\vector(-1,0){16}}
      \put(20,22){\vector(0,1){16}}
      \put(58,40){\line(-1,0){16}}
            \qbezier[10](22,38)(30,30)(38,22)
            \qbezier[10](8,42)(20,47)(35,42)
    }

      \put(260,-5){
      \put(0,20){\circle*{2}}
      \put(20,20){\circle*{2}}
      \put(40,20){\circle*{2}}
      \put(20,5){\circle*{2}}
      \put(40,40){\circle*{2}}
      \put(20,40){\circle*{2}}
      \put(38,20){\vector(-1,0){16}}
      \put(20,7){\line(0,1){11}}
      \put(40,22){\vector(0,1){16}}
      \put(38,40){\vector(-1,0){16}}
      \put(20,22){\vector(0,1){16}}
      \put(2,20){\vector(1,0){16}}
            \qbezier[10](22,38)(30,30)(38,22)
            \qbezier[10](6,22)(17,23)(18,32)
    }

      \put(297,-15){
      \put(60,50){\circle*{2}}
      \put(20,20){\circle*{2}}
      \put(40,20){\circle*{2}}
      \put(60,20){\circle*{2}}
      \put(30,35){\circle*{2}}
      \put(50,35){\circle*{2}}
      \multiput(22,20)(20,0){2}{\vector(1,0){16}}
      \put(38.7,22){\vector(-2,3){7.5}}
      \put(51.7,37){\line(2,3){7.2}}
      \put(48.7,33){\vector(-2,-3){7.5}}
           \qbezier[10](26,21)(40,20)(33,28)
           \qbezier[10](52,21)(40,20)(48,30)
           \qbezier[8](36,30)(40,20)(44,31)
    }
\end{picture}

  \begin{picture}(50,45)

      \put(-1,-5){
      \put(40,5){\circle*{2}}
      \put(20,20){\circle*{2}}
      \put(40,20){\circle*{2}}
      \put(60,20){\circle*{2}}
      \put(40,40){\circle*{2}}
      \put(20,40){\circle*{2}}
      \put(38,20){\vector(-1,0){16}}
      \put(40,7){\vector(0,1){11}}
      \put(40,22){\vector(0,1){16}}
      \put(38,40){\vector(-1,0){16}}
      \put(20,22){\vector(0,1){16}}
      \put(58,20){\vector(-1,0){16}}
            \qbezier[10](22,38)(30,30)(38,22)
            \qbezier[10](42,32)(42,22)(54,22)
    }

      \put(89,-5){
      \put(0,40){\circle*{2}}
      \put(20,20){\circle*{2}}
      \put(40,20){\circle*{2}}
      \put(60,40){\circle*{2}}
      \put(40,40){\circle*{2}}
      \put(20,40){\circle*{2}}
      \put(38,20){\vector(-1,0){16}}
      \put(18,40){\vector(-1,0){16}}
      \put(40,38){\vector(0,-1){16}}
      \put(38,40){\vector(-1,0){16}}
      \put(20,38){\vector(0,-1){16}}
      \put(58,40){\vector(-1,0){16}}
            \qbezier[10](38,38)(30,30)(22,22)
            \qbezier[10](8,42)(20,47)(34,42)
    }

      \put(151,0){
      \put(40,5){\circle*{2}}
      \put(20,20){\circle*{2}}
      \put(40,20){\circle*{2}}
      \put(60,20){\circle*{2}}
      \put(30,35){\circle*{2}}
      \put(50,35){\circle*{2}}
      \multiput(22,20)(20,0){2}{\vector(1,0){16}}
      \put(38.7,22){\vector(-2,3){7.5}}
      \put(40,7){\vector(0,1){11}}
      \put(48.7,33){\vector(-2,-3){7.5}}
           \qbezier[10](25,21)(40,20)(33,28.5)
           \qbezier[10](54,21)(40,20)(47,28.5)
           \qbezier[8](35,30)(40,20)(45,32)
    }


     \put(227,-5){
      \put(60,40){\circle*{2}}
      \put(20,20){\circle*{2}}
      \put(40,20){\circle*{2}}
      \put(60,20){\circle*{2}}
      \put(40,40){\circle*{2}}
      \put(20,40){\circle*{2}}
      \put(38,20){\vector(-1,0){16}}
      \put(42,40){\vector(1,0){16}}
      \put(40,38){\vector(0,-1){16}}
      \put(38,40){\vector(-1,0){16}}
      \put(20,38){\vector(0,-1){16}}
      \put(60,38){\vector(0,-1){16}}
      \put(42,20){\vector(1,0){16}}
            \qbezier[10](22,22)(30,30)(38,38)
            \qbezier[10](58,22)(50,30)(42,38)
    }

     \put(297,-5){
      \put(60,40){\circle*{2}}
      \put(20,20){\circle*{2}}
      \put(40,20){\circle*{2}}
      \put(60,20){\circle*{2}}
      \put(40,40){\circle*{2}}
      \put(20,40){\circle*{2}}
      \put(38,20){\vector(-1,0){16}}
      \put(42,40){\vector(1,0){16}}
      \put(40,22){\vector(0,1){16}}
      \put(60,22){\vector(0,1){16}}
      \put(42,20){\vector(1,0){16}}
      \put(22,38){\vector(1,-1){16}}
            \qbezier[10](42,22)(50,30)(58,38)
            \qbezier[10](28,22)(40,20)(23,34)
            \qbezier[10](26.5,37)(40,20)(38,32)
    }
\end{picture}

  \begin{picture}(50,45)
     \put(4,-5){
      \put(60,40){\circle*{2}}
      \put(80,20){\circle*{2}}
      \put(40,20){\circle*{2}}
      \put(60,20){\circle*{2}}
      \put(40,40){\circle*{2}}
      \put(20,40){\circle*{2}}
      \put(22,40){\vector(1,0){16}}
      \put(42,40){\vector(1,0){16}}
      \put(40,22){\vector(0,1){16}}
      \put(60,22){\vector(0,1){16}}
      \put(42,20){\vector(1,0){16}}
      \put(78,20){\vector(-1,0){16}}
            \qbezier[10](42,22)(50,30)(58,38)
            \qbezier[10](26,42)(40,47)(52,42)
            \qbezier[10](62,32)(64,23)(74,22)
    }

     \put(80,-20){
      \put(60,40){\circle*{2}}
      \put(80,40){\circle*{2}}
      \put(40,20){\circle*{2}}
      \put(60,20){\circle*{2}}
      \put(40,40){\circle*{2}}
      \put(60,55){\circle*{2}}
      \put(60,42){\vector(0,1){11}}
      \put(42,40){\vector(1,0){16}}
      \put(40,22){\vector(0,1){16}}
      \put(60,22){\vector(0,1){16}}
      \put(42,20){\vector(1,0){16}}
      \put(62,40){\vector(1,0){16}}
            \qbezier[10](46,42)(58,42)(58,48)
            \qbezier[10](42,22)(50,30)(58,38)
            \qbezier[10](62,24)(64,38)(72,38)
    }

      \put(190,-5){
      \put(60,40){\circle*{2}}
      \put(80,40){\circle*{2}}
      \put(40,20){\circle*{2}}
      \put(60,20){\circle*{2}}
      \put(40,40){\circle*{2}}
      \put(20,40){\circle*{2}}
      \put(22,40){\vector(1,0){16}}
      \put(42,40){\vector(1,0){16}}
      \put(40,22){\vector(0,1){16}}
      \put(60,22){\vector(0,1){16}}
      \put(42,20){\vector(1,0){16}}
      \put(62,40){\vector(1,0){16}}
            \qbezier[10](42,22)(50,30)(58,38)
            \qbezier[10](26,42)(40,47)(52,42)
            \qbezier[10](62,26)(64,37)(72,38)
    }

      \put(297,-5){
      \put(60,40){\circle*{2}}
      \put(20,40){\circle*{2}}
      \put(40,20){\circle*{2}}
      \put(60,20){\circle*{2}}
      \put(40,40){\circle*{2}}
      \put(20,20){\circle*{2}}
      \put(38,20){\vector(-1,0){16}}
      \put(22,38){\vector(1,-1){16}}
      \put(40,22){\vector(0,1){16}}
      \put(58,38){\vector(-1,-1){16}}
      \put(42,20){\vector(1,0){16}}
            \qbezier[10](42,32)(40,20)(52,35)
            \qbezier[10](38,32)(40,20)(28,35)
            \qbezier[10](52,22)(40,20)(57,34)
            \qbezier[10](28,22)(40,20)(24,33)
    }
\end{picture}

  \begin{picture}(50,35)

     \put(3,-15){
      \put(60,40){\circle*{2}}
      \put(20,20){\circle*{2}}
      \put(40,20){\circle*{2}}
      \put(60,20){\circle*{2}}
      \put(40,40){\circle*{2}}
      \put(20,40){\circle*{2}}
      \put(22,20){\vector(1,0){16}}
      \put(42,40){\vector(1,0){16}}
      \put(40,38){\vector(0,-1){16}}
      \put(22,40){\vector(1,0){16}}
      \put(20,38){\vector(0,-1){16}}
      \put(60,38){\vector(0,-1){16}}
      \put(42,20){\vector(1,0){16}}
            \qbezier[10](22,38)(30,30)(38,22)
            \qbezier[10](58,22)(50,30)(42,38)
    }

      \put(85,-15){
      \put(60,40){\circle*{2}}
      \put(80,40){\circle*{2}}
      \put(40,20){\circle*{2}}
      \put(60,20){\circle*{2}}
      \put(40,40){\circle*{2}}
      \put(20,40){\circle*{2}}
      \put(38,40){\vector(-1,0){16}}
      \put(42,40){\vector(1,0){16}}
      \put(40,22){\vector(0,1){16}}
      \put(60,22){\vector(0,1){16}}
      \put(42,20){\vector(1,0){16}}
      \put(62,40){\vector(1,0){16}}
            \qbezier[10](42,22)(50,30)(58,38)
            \qbezier[10](28,38)(38,37)(38,26)
            \qbezier[10](62,26)(64,37)(72,38)
    }

     \put(190,-15){
      \put(60,40){\circle*{2}}
      \put(80,20){\circle*{2}}
      \put(40,20){\circle*{2}}
      \put(60,20){\circle*{2}}
      \put(40,40){\circle*{2}}
      \put(20,40){\circle*{2}}
      \put(38,40){\vector(-1,0){16}}
      \put(42,40){\vector(1,0){16}}
      \put(40,22){\vector(0,1){16}}
      \put(60,22){\vector(0,1){16}}
      \put(42,20){\vector(1,0){16}}
      \put(78,20){\vector(-1,0){16}}
            \qbezier[10](42,22)(50,30)(58,38)
            \qbezier[10](28,38)(38,37)(38,26)
            \qbezier[10](62,32)(64,23)(74,22)
    }

     \put(297,-15){
      \put(60,40){\circle*{2}}
      \put(20,20){\circle*{2}}
      \put(40,20){\circle*{2}}
      \put(60,20){\circle*{2}}
      \put(40,40){\circle*{2}}
      \put(20,40){\circle*{2}}
      \put(38,20){\vector(-1,0){16}}
      \put(58,40){\vector(-1,0){16}}
      \put(40,22){\vector(0,1){16}}
      \put(60,22){\vector(0,1){16}}
      \put(58,20){\vector(-1,0){16}}
      \put(22,38){\vector(1,-1){16}}
            \qbezier[10](58,22)(50,30)(42,38)
            \qbezier[10](28,22)(40,20)(23,34)
            \qbezier[10](26.5,37)(40,20)(38,32)
    }
\end{picture}

  \begin{picture}(50,45)

      \put(3,-5){
      \put(60,40){\circle*{2}}
      \put(20,20){\circle*{2}}
      \put(40,20){\circle*{2}}
      \put(60,20){\circle*{2}}
      \put(40,40){\circle*{2}}
      \put(20,40){\circle*{2}}
      \put(22,20){\vector(1,0){16}}
      \put(42,40){\vector(1,0){16}}
      \put(38,40){\vector(-1,0){16}}
      \put(20,38){\vector(0,-1){16}}
       \put(60,38){\vector(0,-1){16}}
      \put(58,20){\vector(-1,0){16}}
             \qbezier[8](40,22)(40,30)(40,38)

    }

      \put(92,0){
      \put(80,20){\circle*{2}}
      \put(20,20){\circle*{2}}
      \put(40,20){\circle*{2}}
      \put(60,20){\circle*{2}}
      \put(40,35){\circle*{2}}
      \put(0,20){\circle*{2}}
      \put(22,20){\vector(1,0){16}}
      \put(78,20){\vector(-1,0){16}}
      \put(40,22){\vector(0,1){11}}
      \put(2,20){\vector(1,0){16}}
      \put(58,20){\vector(-1,0){16}}
             \qbezier[15](6,22)(35,22)(38,28)
             \qbezier[15](74,22)(45,22)(42,28)
    }

      \put(200,-15){
      \put(40,50){\circle*{2}}
      \put(20,20){\circle*{2}}
      \put(40,20){\circle*{2}}
      \put(60,20){\circle*{2}}
      \put(40,35){\circle*{2}}
      \put(0,20){\circle*{2}}
      \put(22,20){\vector(1,0){16}}
      \put(40,37){\vector(0,1){11}}
      \put(40,22){\vector(0,1){11}}
      \put(2,20){\vector(1,0){16}}
      \put(58,20){\vector(-1,0){16}}
             \qbezier[15](6,22)(35,22)(38,28)
             \qbezier[15](54,22)(42,25)(42,43)
    }

      \put(295,0){
      \put(20,5){\circle*{2}}
      \put(20,20){\circle*{2}}
      \put(40,20){\circle*{2}}
      \put(60,20){\circle*{2}}
      \put(40,35){\circle*{2}}
      \put(0,20){\circle*{2}}
      \put(22,20){\vector(1,0){16}}
      \put(20,7){\vector(0,1){11}}
      \put(40,22){\vector(0,1){11}}
      \put(2,20){\vector(1,0){16}}
      \put(42,20){\vector(1,0){16}}
             \qbezier[15](6,22)(35,22)(38,28)
             \qbezier[15](22,10)(25,18)(52,18)
    }
\end{picture}

  \begin{picture}(50,50)

      \put(7,-15){
      \put(0,20){\circle*{2}}
      \put(20,20){\circle*{2}}
      \put(40,20){\circle*{2}}
      \put(30,55){\circle*{2}}
      \put(40,40){\circle*{2}}
      \put(20,40){\circle*{2}}
      \put(22,20){\vector(1,0){16}}
      \put(18,20){\line(-1,0){16}}
      \put(40,22){\vector(0,1){16}}
      \put(20,38){\vector(0,-1){16}}
      \put(21,42){\vector(2,3){7.7}}
      \put(31.7,53.5){\vector(2,-3){7.7}}
            \qbezier[8](22,40)(30,40)(38,40)
    }

      \put(90,5){
      \put(-20,20){\circle*{2}}
      \put(20,20){\circle*{2}}
      \put(40,20){\circle*{2}}
      \put(60,20){\circle*{2}}
      \put(40,35){\circle*{2}}
      \put(0,20){\circle*{2}}
      \put(22,20){\vector(1,0){16}}
      \put(-18,20){\line(1,0){16}}
      \put(40,22){\vector(0,1){11}}
      \put(2,20){\vector(1,0){16}}
      \put(58,20){\vector(-1,0){16}}
             \qbezier[15](6,22)(35,22)(38,28)
             \qbezier[10](42,28)(43,23)(54,22)
    }

      \put(175,5){
      \put(20,5){\circle*{2}}
      \put(20,20){\circle*{2}}
      \put(40,20){\circle*{2}}
      \put(60,20){\circle*{2}}
      \put(0,20){\circle*{2}}
      \put(80,20){\circle*{2}}
      \put(22,20){\vector(1,0){16}}
      \put(20,7){\vector(0,1){11}}
      \put(2,20){\vector(1,0){16}}
      \put(62,20){\line(1,0){16}}
      \put(42,20){\vector(1,0){16}}
             \qbezier[12](6,22)(15,29)(32,22)
             \qbezier[12](22,10)(25,18)(52,18)
    }

      \put(275,5){
      \put(40,35){\circle*{2}}
      \put(20,20){\circle*{2}}
      \put(40,20){\circle*{2}}
      \put(60,20){\circle*{2}}
      \put(0,20){\circle*{2}}
      \put(80,20){\circle*{2}}
      \put(22,20){\vector(1,0){16}}
      \put(2,20){\vector(1,0){16}}
      \put(40,22){\vector(0,1){11}}
      \put(62,20){\line(1,0){16}}
      \put(42,20){\vector(1,0){16}}
             \qbezier[15](6,18)(25,11)(52,18)
             \qbezier[10](25,22)(38,22)(38,28)
    }
\end{picture}

  \begin{picture}(50,45)

      \put(6,0){
      \put(20,5){\circle*{2}}
      \put(20,20){\circle*{2}}
      \put(40,20){\circle*{2}}
      \put(60,20){\circle*{2}}
      \put(40,35){\circle*{2}}
      \put(0,20){\circle*{2}}
      \put(18,20){\vector(-1,0){16}}
      \put(20,7){\line(0,1){11}}
      \put(40,33){\vector(0,-1){11}}
      \put(38,20){\vector(-1,0){16}}
      \put(42,20){\vector(1,0){16}}
             \qbezier[15](8,22)(35,22)(38,30)
             \qbezier[10](42,30)(43,23)(52,22)
    }


       \put(87,-20){
      \put(0,40){\circle*{2}}
      \put(20,20){\circle*{2}}
      \put(40,20){\circle*{2}}
      \put(30,55){\circle*{2}}
      \put(40,40){\circle*{2}}
      \put(20,40){\circle*{2}}
      \put(22,20){\vector(1,0){16}}
      \put(18,40){\line(-1,0){16}}
      \put(40,22){\vector(0,1){16}}
      \put(20,38){\vector(0,-1){16}}
      \put(21,42){\vector(2,3){7.7}}
      \put(31.7,53.5){\vector(2,-3){7.7}}
            \qbezier[8](22,40)(30,40)(38,40)
    }

      \put(145,-15){
      \put(40,50){\circle*{2}}
      \put(20,20){\circle*{2}}
      \put(40,20){\circle*{2}}
      \put(60,20){\circle*{2}}
      \put(40,35){\circle*{2}}
      \put(0,20){\circle*{2}}
      \put(2,20){\vector(1,0){16}}
      \put(40,37){\line(0,1){11}}
      \put(40,22){\vector(0,1){11}}
      \put(22,20){\vector(1,0){16}}
      \put(58,20){\vector(-1,0){16}}
             \qbezier[15](6,22)(35,22)(38,28)
             \qbezier[10](42,28)(42,22)(54,22)
    }

      \put(220,0){
      \put(20,5){\circle*{2}}
      \put(20,20){\circle*{2}}
      \put(40,20){\circle*{2}}
      \put(60,20){\circle*{2}}
      \put(40,35){\circle*{2}}
      \put(0,20){\circle*{2}}
      \put(2,20){\vector(1,0){16}}
      \put(20,7){\vector(0,1){11}}
      \put(40,22){\vector(0,1){11}}
      \put(22,20){\vector(1,0){16}}
      \put(58,20){\line(-1,0){16}}
             \qbezier[15](6,22)(35,22)(38,28)
             \qbezier[10](22,10)(22,18)(32,18)
    }

      \put(295,0){
      \put(40,5){\circle*{2}}
      \put(20,20){\circle*{2}}
      \put(40,20){\circle*{2}}
      \put(60,20){\circle*{2}}
      \put(40,35){\circle*{2}}
      \put(0,20){\circle*{2}}
      \put(2,20){\vector(1,0){16}}
      \put(40,18){\vector(0,-1){11}}
      \put(40,22){\vector(0,1){11}}
      \put(22,20){\vector(1,0){16}}
      \put(58,20){\line(-1,0){16}}
             \qbezier[15](6,22)(35,22)(38,28)
             \qbezier[10](25,18)(38,18)(38,12)
    }

\end{picture}

  \begin{picture}(50,50)

      \put(10,0){
      \put(40,5){\circle*{2}}
      \put(20,20){\circle*{2}}
      \put(40,20){\circle*{2}}
      \put(60,20){\circle*{2}}
      \put(40,35){\circle*{2}}
      \put(0,20){\circle*{2}}
      \put(18,20){\vector(-1,0){16}}
      \put(40,7){\line(0,1){11}}
      \put(40,33){\vector(0,-1){11}}
      \put(38,20){\vector(-1,0){16}}
      \put(42,20){\vector(1,0){16}}
             \qbezier[15](8,22)(35,22)(38,30)
             \qbezier[10](42,30)(43,23)(52,22)
    }


      \put(110,-5){
      \put(60,40){\circle*{2}}
      \put(0,40){\circle*{2}}
      \put(40,20){\circle*{2}}
      \put(60,20){\circle*{2}}
      \put(40,40){\circle*{2}}
      \put(20,40){\circle*{2}}
      \put(18,40){\vector(-1,0){16}}
      \put(58,40){\vector(-1,0){16}}
      \put(22,40){\vector(1,0){16}}
      \put(38,22){\vector(-1,1){16}}
      \put(60,22){\vector(0,1){16}}
      \put(42,20){\vector(1,0){16}}
             \qbezier[8](40,22)(40,30)(40,38)
             \qbezier[12](7,38)(20,37)(35,22)
    }

      \put(196,-5){
      \put(60,40){\circle*{2}}
      \put(40,5){\circle*{2}}
      \put(40,20){\circle*{2}}
      \put(60,20){\circle*{2}}
      \put(40,40){\circle*{2}}
      \put(20,40){\circle*{2}}
      \put(22,40){\vector(1,0){16}}
      \put(58,40){\vector(-1,0){16}}
      \put(40,7){\vector(0,1){11}}
      \put(38,22){\vector(-1,1){16}}
      \put(60,22){\vector(0,1){16}}
      \put(42,20){\vector(1,0){16}}
             \qbezier[8](40,22)(40,30)(40,38)
             \qbezier[12](38,10)(37,20)(25,32)
    }

         \put(295,-5){
      \put(60,40){\circle*{2}}
      \put(20,20){\circle*{2}}
      \put(40,20){\circle*{2}}
      \put(60,20){\circle*{2}}
      \put(40,40){\circle*{2}}
      \put(0,20){\circle*{2}}
      \put(22,20){\vector(1,0){16}}
      \put(58,40){\vector(-1,0){16}}
      \put(40,22){\vector(0,1){16}}
      \put(60,22){\vector(0,1){16}}
      \put(58,20){\vector(-1,0){16}}
      \put(2,20){\vector(1,0){16}}
            \qbezier[10](58,22)(50,30)(42,38)
            \qbezier[15](6,22)(35,22)(38,32)
    }
\end{picture}

  \begin{picture}(50,50)
         \put(10,0){
      \put(60,40){\circle*{2}}
      \put(20,20){\circle*{2}}
      \put(40,20){\circle*{2}}
      \put(60,20){\circle*{2}}
      \put(40,40){\circle*{2}}
      \put(0,20){\circle*{2}}
      \put(38,20){\vector(-1,0){16}}
      \put(58,40){\vector(-1,0){16}}
      \put(40,38){\vector(0,-1){16}}
      \put(60,38){\vector(0,-1){16}}
      \put(58,20){\vector(-1,0){16}}
      \put(18,20){\vector(-1,0){16}}
            \qbezier[10](42,22)(50,30)(58,38)
            \qbezier[15](8,22)(35,22)(38,34)
    }

         \put(80,-15){
      \put(60,40){\circle*{2}}
      \put(20,20){\circle*{2}}
      \put(40,20){\circle*{2}}
      \put(60,20){\circle*{2}}
      \put(40,40){\circle*{2}}
      \put(40,55){\circle*{2}}
      \put(22,20){\vector(1,0){16}}
      \put(58,40){\vector(-1,0){16}}
      \put(40,22){\vector(0,1){16}}
      \put(60,22){\vector(0,1){16}}
      \put(58,20){\vector(-1,0){16}}
      \put(40,42){\vector(0,1){11}}
            \qbezier[10](58,22)(50,30)(42,38)
            \qbezier[15](24,22)(38,25)(38,47)
    }

     \put(180,0){
      \put(80,20){\circle*{2}}
      \put(20,20){\circle*{2}}
      \put(40,20){\circle*{2}}
      \put(60,20){\circle*{2}}
      \put(40,40){\circle*{2}}
      \put(20,40){\circle*{2}}
      \put(22,20){\vector(1,0){16}}
      \put(40,38){\vector(0,-1){16}}
      \put(62,20){\vector(1,0){16}}
      \put(42,20){\vector(1,0){16}}
      \put(38,22){\vector(-1,1){16}}
            \qbezier[15](42,36)(45,22)(72,22)
            \qbezier[10](25,22)(40,20)(25,33)
            \qbezier[10](27,36)(40,20)(38,36)
    }

      \put(295,-15){
      \put(40,55){\circle*{2}}
      \put(20,20){\circle*{2}}
      \put(40,20){\circle*{2}}
      \put(60,20){\circle*{2}}
      \put(40,40){\circle*{2}}
      \put(20,40){\circle*{2}}
      \put(22,20){\vector(1,0){16}}
      \put(40,53){\vector(0,-1){11}}
      \put(40,38){\vector(0,-1){16}}
      \put(42,20){\vector(1,0){16}}
      \put(38,22){\vector(-1,1){16}}
            \qbezier[15](42,50)(42,25)(52,22)
            \qbezier[10](25,22)(40,20)(25,33)
            \qbezier[10](27,36)(40,20)(38,36)
    }

  \end{picture}
\end{thm}\vspace{.1in}

The characterization of tilted algebras of type $A_n$ in terms of
their quiver with relations has been done by I. Assem in \cite{a}.
F. Huard, in \cite{huard}, has characterized the tilted algebras
$kQ/I$ where the underlying graph of $Q$ is $D_n$. An alternative
proof of this result can be made using the method described in this
work. A complete classification  of tilted algebras of type $E_p,
p=7,8$, was done in \cite{bft2}.


\subsection{Cluster tilted algebras of finite representation type}

The relationship between tilted algebras  and cluster-tilted
algebras was given in \cite{abs}. For this, the authors introduced
the concept of relation extension $R(B)$ of an algebra $B$ with
$\mathrm{gl.dim.}B \leq 2$ . Recall that the algebra $R(B) = B
\propto \mathrm{Ext}^2_B(DB,B)$, where $DB = \mathrm{Hom}_k(B, k)$.
In case the quiver of $B$ has no oriented cycles,  a construction of
the quiver of the relation extension algebra is given in \cite{abs}
. The main result in \cite{abs} characterizes the cluster tilted
algebras as the relation-extension of
 tilted algebras.

On the other hand, A. Buan, R. Marsh and I.Reiten in \cite{bmr2} proved
that any cluster tilted algebra of finite representation type is
uniquely determined by its ordinary quiver (up to isomorphism).

Using Theorem \ref{thm:te6} and the two previous
results, we deduce the following Theorem.

\begin{thm}\label{thm:e6}
Let $C$ be an algebra. Then $C$ is cluster tilted of Dynkin type
$E_6$ if and only if
\begin{enumerate}
    \item[i.] $C$ is a hereditary algebra of Dynkin type $E_6$,
    or
    \item[ii.] the quiver $Q_C$ or $Q_{C^{op}}$ is isomorphic to one of the following:
\end{enumerate}
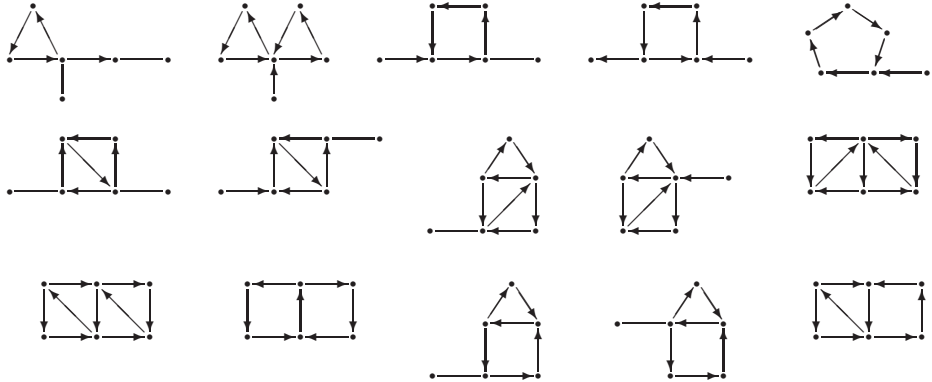
\vspace{.05in}

  \begin{picture}(50,40)

      \put(7,-5){
      \put(0,20){\circle*{2}}
      \put(9,40){\circle*{2}}
      \put(20,20){\circle*{2}}
      \put(60,20){\circle*{2}}
      \put(40,20){\circle*{2}}
      \put(20,5){\circle*{2}}
      \put(8,38){\vector(-1,-2){8}}
      \put(18,22){\vector(-1,2){8}}
      \multiput(2,20)(20,0){2}{\vector(1,0){16}}
      \put(20,7){\line(0,2){11}}
      \put(42,20){\line(1,0){16}}
      }

      \put(87,-5){
      \put(0,20){\circle*{2}}
      \put(9,40){\circle*{2}}
      \put(20,20){\circle*{2}}
      \put(30,40){\circle*{2}}
      \put(40,20){\circle*{2}}
      \put(20,5){\circle*{2}}
      \put(8,38){\vector(-1,-2){8}}
      \put(18,22){\vector(-1,2){8}}
      \put(38,22){\vector(-1,2){8}}
      \multiput(2,20)(20,0){2}{\vector(1,0){16}}
      \put(28,38){\vector(-1,-2){8}}
      \put(20,7){\vector(0,2){11}}
    }

      \put(147,-5){
      \put(0,20){\circle*{2}}
      \put(20,20){\circle*{2}}
      \put(40,20){\circle*{2}}
      \put(60,20){\circle*{2}}
      \put(40,40){\circle*{2}}
      \put(20,40){\circle*{2}}
      \multiput(2,20)(20,0){2}{\vector(1,0){16}}
      \put(42,20){\line(1,0){16}}
      \put(40,22){\vector(0,1){16}}
      \put(38,40){\vector(-1,0){16}}
      \put(20,38){\vector(0,-1){16}}
    }

      \put(227,-5){
      \put(0,20){\circle*{2}}
      \put(20,20){\circle*{2}}
      \put(40,20){\circle*{2}}
      \put(60,20){\circle*{2}}
      \put(40,40){\circle*{2}}
      \put(20,40){\circle*{2}}
      \put(18,20){\vector(-1,0){16}}
      \put(58,20){\vector(-1,0){16}}
      \put(22,20){\vector(1,0){16}}
      \put(40,22){\vector(0,1){16}}
      \put(38,40){\vector(-1,0){16}}
      \put(20,38){\vector(0,-1){16}}
    }

      \put(314,-10){
      \put(-5,35){\circle*{2}}
      \put(0,20){\circle*{2}}
      \put(10,45){\circle*{2}}
      \put(25,35){\circle*{2}}
      \put(20,20){\circle*{2}}
      \put(40,20){\circle*{2}}
      \put(-0.5,22){\vector(-1,3){3.5}}
      \put(-3,36){\vector(3,2){11}}
      \put(12,44){\vector(3,-2){11}}
      \put(24,33){\vector(-1,-3){3.5}}
      \put(18,20){\vector(-1,0){16}}
      \put(38,20){\vector(-1,0){16}}
    }
\end{picture}

  \begin{picture}(50,50)
      \put(7,0){
      \put(0,20){\circle*{2}}
      \put(20,20){\circle*{2}}
      \put(40,20){\circle*{2}}
      \put(60,20){\circle*{2}}
      \put(40,40){\circle*{2}}
      \put(20,40){\circle*{2}}
      \put(38,20){\vector(-1,0){16}}
      \put(42,20){\line(1,0){16}}
      \put(18,20){\line(-1,0){16}}
      \put(40,22){\vector(0,1){16}}
      \put(38,40){\vector(-1,0){16}}
      \put(20,22){\vector(0,1){16}}
      \put(22,38){\vector(1,-1){16}}
    }

      \put(87,0){
      \put(0,20){\circle*{2}}
      \put(20,20){\circle*{2}}
      \put(40,20){\circle*{2}}
      \put(60,40){\circle*{2}}
      \put(40,40){\circle*{2}}
      \put(20,40){\circle*{2}}
      \put(38,20){\vector(-1,0){16}}
      \put(42,40){\line(1,0){16}}
      \put(2,20){\vector(1,0){16}}
      \put(40,22){\vector(0,1){16}}
      \put(38,40){\vector(-1,0){16}}
      \put(20,22){\vector(0,1){16}}
      \put(22,38){\vector(1,-1){16}}
    }

      \put(166,-15){
      \put(0,20){\circle*{2}}
      \put(20,20){\circle*{2}}
      \put(40,20){\circle*{2}}
      \put(30,55){\circle*{2}}
      \put(40,40){\circle*{2}}
      \put(20,40){\circle*{2}}
      \put(38,20){\vector(-1,0){16}}
      \put(18,20){\line(-1,0){16}}
      \put(40,38){\vector(0,-1){16}}
      \put(38,40){\vector(-1,0){16}}
      \put(20,38){\vector(0,-1){16}}
      \put(22,22){\vector(1,1){16}}
      \put(21,42){\vector(2,3){7.7}}
      \put(31.7,53.5){\vector(2,-3){7.7}}
    }

      \put(219,-15){
      \put(60,40){\circle*{2}}
      \put(20,20){\circle*{2}}
      \put(40,20){\circle*{2}}
      \put(30,55){\circle*{2}}
      \put(40,40){\circle*{2}}
      \put(20,40){\circle*{2}}
      \put(38,20){\vector(-1,0){16}}
      \put(58,40){\vector(-1,0){16}}
      \put(40,38){\vector(0,-1){16}}
      \put(38,40){\vector(-1,0){16}}
      \put(20,38){\vector(0,-1){16}}
      \put(22,22){\vector(1,1){16}}
      \put(21,42){\vector(2,3){7.7}}
      \put(31.7,53.5){\vector(2,-3){7.7}}
    }

     \put(290,0){
      \put(60,40){\circle*{2}}
      \put(20,20){\circle*{2}}
      \put(40,20){\circle*{2}}
      \put(60,20){\circle*{2}}
      \put(40,40){\circle*{2}}
      \put(20,40){\circle*{2}}
      \put(38,20){\vector(-1,0){16}}
      \put(42,40){\vector(1,0){16}}
      \put(40,38){\vector(0,-1){16}}
      \put(38,40){\vector(-1,0){16}}
      \put(20,38){\vector(0,-1){16}}
      \put(22,22){\vector(1,1){16}}
      \put(60,38){\vector(0,-1){16}}
      \put(42,20){\vector(1,0){16}}
      \put(58,22){\vector(-1,1){16}}
    }
\end{picture}

  \begin{picture}(50,50)

     \put(0,0){
      \put(60,40){\circle*{2}}
      \put(20,20){\circle*{2}}
      \put(40,20){\circle*{2}}
      \put(60,20){\circle*{2}}
      \put(40,40){\circle*{2}}
      \put(20,40){\circle*{2}}
      \put(22,20){\vector(1,0){16}}
      \put(42,40){\vector(1,0){16}}
      \put(40,38){\vector(0,-1){16}}
      \put(22,40){\vector(1,0){16}}
      \put(20,38){\vector(0,-1){16}}
      \put(38,22){\vector(-1,1){16}}
      \put(60,38){\vector(0,-1){16}}
      \put(42,20){\vector(1,0){16}}
      \put(58,22){\vector(-1,1){16}}
    }

     \put(77,0){
      \put(60,40){\circle*{2}}
      \put(20,20){\circle*{2}}
      \put(40,20){\circle*{2}}
      \put(60,20){\circle*{2}}
      \put(40,40){\circle*{2}}
      \put(20,40){\circle*{2}}
      \put(22,20){\vector(1,0){16}}
      \put(42,40){\vector(1,0){16}}
      \put(40,22){\vector(0,1){16}}
      \put(38,40){\vector(-1,0){16}}
      \put(20,38){\vector(0,-1){16}}
       \put(60,38){\vector(0,-1){16}}
      \put(58,20){\vector(-1,0){16}}
    }

      \put(167,-15){
      \put(0,20){\circle*{2}}
      \put(20,20){\circle*{2}}
      \put(40,20){\circle*{2}}
      \put(30,55){\circle*{2}}
      \put(40,40){\circle*{2}}
      \put(20,40){\circle*{2}}
      \put(22,20){\vector(1,0){16}}
      \put(18,20){\line(-1,0){16}}
      \put(40,22){\vector(0,1){16}}
      \put(38,40){\vector(-1,0){16}}
      \put(20,38){\vector(0,-1){16}}
      \put(21,42){\vector(2,3){7.7}}
      \put(31.7,53.5){\vector(2,-3){7.7}}
    }

      \put(237,-15){
      \put(0,40){\circle*{2}}
      \put(20,20){\circle*{2}}
      \put(40,20){\circle*{2}}
      \put(30,55){\circle*{2}}
      \put(40,40){\circle*{2}}
      \put(20,40){\circle*{2}}
      \put(22,20){\vector(1,0){16}}
      \put(18,40){\line(-1,0){16}}
      \put(40,22){\vector(0,1){16}}
      \put(38,40){\vector(-1,0){16}}
      \put(20,38){\vector(0,-1){16}}
      \put(21,42){\vector(2,3){7.7}}
      \put(31.7,53.5){\vector(2,-3){7.7}}
    }

     \put(292,0){
      \put(60,40){\circle*{2}}
      \put(20,20){\circle*{2}}
      \put(40,20){\circle*{2}}
      \put(60,20){\circle*{2}}
      \put(40,40){\circle*{2}}
      \put(20,40){\circle*{2}}
      \put(22,20){\vector(1,0){16}}
      \put(58,40){\vector(-1,0){16}}
      \put(40,38){\vector(0,-1){16}}
      \put(22,40){\vector(1,0){16}}
      \put(20,38){\vector(0,-1){16}}
      \put(38,22){\vector(-1,1){16}}
      \put(60,22){\vector(0,1){16}}
      \put(42,20){\vector(1,0){16}}
    }

  \end{picture}

\end{thm}

It is interesting to point out that, the classification of the
cluster tilted algebras of the remaining Dynkin types, can be done
in a similar way. Since the cluster tilted algebras of type $A_n$
and $D_n$ were studied in \cite{ccs} and \cite{s} respectively, by
geometric methods, and in \cite{fst} was shown that it is not
possible to obtain the $E_p$ type using geometric methods, we are
mainly interested in the cluster tilted algebras of type $E_p,\:
p=7,8 $. For classification of these two classes of algebras see
\cite{bft2}.



\begin{thebibliography}{AST}

\bibitem{a} I. Assem. {\em Tilted
algebras of type $A_n$}. Comm. Algebra 19, (1982), 2121-2131.

\bibitem{abs} I. Assem, T. Br$\mathrm{\ddot{u}}$stle, R. Schiffler. {\em Cluster-tilted
algebras as trivial extensions}. Bull. Lond. Math. Soc. 40 (2008),
no. 1, 151-162.

\bibitem{ar} I. Assem, M. J. Redondo. {\em The first Hochschild cohomology group of a
schurian cluster-tilted algebra}. Manuscripta Math. 128, (2009),
no. 3, 373-388.


\bibitem{ars}  M. Auslander, I. Reiten, S. Smal\o. {\em
Representation Theory of Artin Algebras}, Cambridge Studies in
Advanced Mathematics,  vol 36, Cambridge University Press (1995).

\bibitem{bg}  K. Bongartz, P. Gabriel. {\em Covering Spaces in
Representation-Theory}. Invent. Mathematicae 65, (1982), 331-378.

\bibitem{bft} N. Bordino, E.  Fern\'andez, S. Trepode. {\em A criterion for global dimension two for strongly
simply connected schurian algebras}.    arXiv:1110.6160 [math.RA].

\bibitem{bft2} N. Bordino, E.  Fern\'andez, S. Trepode. {\em Quivers with relations of
tilted algebras and cluster tilted algebras of type $E_p$}. In preparation.


\bibitem{bmr2} A. Buan, R. Marsh, I. Reiten. {\em Cluster-tilted algebras of finite representation
type}. J. Algebra 306, (2006), no. 2, 412-431.

\bibitem{ccs} P. Caldero, F. Chapoton, R. Schiffler. {\em Quivers with relations arising from clusters
($A_n$ case)}. Trans. Amer. Math. Soc. 358, (2006), no. 3,
1347-1364.

\bibitem{cpt1} F. Coelho, J. A. de la Pe\~na, S. Trepode. {\em Non-quasitilted til
critical algebras}. XVI Coloquio Latinoamericano de \'Algebra,
Proceedings 2007.

\bibitem{cpt2} F. Coelho, J. A. de la Pe\~na, S. Trepode. {\em Minimal non
tilted algebras}. Coloquium Mathematicum 111, (2008), 71-81.

\bibitem{cpt3} F. Coelho, J. A. de la Pe\~na, S. Trepode.
{\em Tilt-critical algebras of tame type}.  Comm. Algebra. Volume
36, (2008), no. 9,  3574-3584.

\bibitem{f} E. Fern\'andez. {\em Extensiones triviales y \'algebras inclinadas iteradas}. Ph.D. Thesis, Universidad
Nacional del Sur, Argentina (1999).

\bibitem{fp} E. Fern\'andez, M. I. Platzeck. {\em Presentations of trivial extensions
of finite dimensional algebras and a theorem of Sheila Brenner}.
J. Algebra 249, (2002), 326-344.

\bibitem{fst} S. Fomin, M. Shapiro, D. Thurston. {\em Cluster algebras and
triangulated surfaces. Part I: Cluster complexes}. Acta Math. 201,
(2008), 83 - 146.

\bibitem{hr} D. Happel, C. Ringel. {\em Tilted
algebras}. Trans. Amer. Math. Soc. 274, (1982), no. 2, 399-443.

\bibitem{hrs} D. Happel, I. Reiten, S. O. Smal\o.
{\em Tilting in Abelian Categories and Quasitilted algebras}.
Memoirs Amer. Math. Soc. 575, (1996).

\bibitem{huard} F. Huard. {\em Tilted algebras having underlying graph $D_n$}.
Annales des Sciences Math\'ematiques du Qu\'ebec (2001), 25 (1)
(2001) 39-46.

\bibitem{h} F. Huard. {\em Tilted gentle algebras}.
Comm. Algebra 26, (1998), no. 1, 63-72.

\bibitem{hl1} F. Huard, S. Liu. {\em Tilted special biserial algebras}.
J. Algebra (1999), 679-700.

\bibitem{hl2} F. Huard, S. Liu. {\em Tilted string algebras}.
J. Pure Appl. Algebra 153, (2000), 151-164.


\bibitem{s} R. Schiffler. {\em A geometric model for cluster categories of type $D_n$}. J. Algebraic Combin. 27, (2008),
no. 1, 1-21.

\bibitem{Sko} A.~Skowro\'nski. {\em Simply connected algebras and Hochschild
  cohomologies.\/}  Representations of algebras (Ottawa, ON, 1992),
431--447, CMS Conf. Proc., 14, Amer. Math. Soc., Providence, RI,
1993.

\end{thebibliography}
\end{document}